  \documentclass[a4paper]{svjour3}       
  \usepackage[margin=1in]{geometry}
  \smartqed  
  \usepackage[utf8]{inputenc}
  
  \usepackage{float, enumerate, graphicx, amssymb, amsmath , color, longtable, multirow, lineno,url}
  \usepackage{graphicx}
  \usepackage{amsmath}
  \usepackage[active]{srcltx}
  \usepackage{pifont}
  \usepackage{dsfont}
  \usepackage{paralist} 
  \usepackage{amssymb}
  \usepackage{authblk}
  \usepackage{epsfig}
  \usepackage{epstopdf}
  \usepackage{subfigure}
  \usepackage{longtable,booktabs}
  \usepackage{mathrsfs}
  \usepackage{lineno}
  
  \usepackage{cite}
  \usepackage{algorithm}
  \usepackage{algorithmic}
  \usepackage[colorlinks,linkcolor=blue,citecolor=blue, bookmarksnumbered=true, bookmarksopen=true, colorlinks ]{hyperref}
  \usepackage{epstopdf}
  \numberwithin{theorem}{section}
  \numberwithin{lemma}{section}
  \numberwithin{corollary}{section}
  \numberwithin{proposition}{section}
  \numberwithin{definition}{section}
  \numberwithin{example}{section}
  \numberwithin{assumption}{section}
  \usepackage{setspace}
  
  \journalname{   }
  
\begin{document}

\title{\textsc{Two  efficient gradient methods with approximately optimal stepsizes based on   regularization models for  unconstrained optimization }  }
\author{ Zexian Liu$ ^{1} $ \and Wangli Chu$ ^{2} $    \and  Hongwei Liu$ ^{2, *} $    }
\institute{
Zexian Liu,   e-mail: liuzexian2008@163.com
\at  School of Mathematics and Statistics, Guizhou University, Guiyang, 550025, China
\and
 Wangli Chu,  e-mail:  chuwangli123@163.com; Hongwei Liu(\ding{41}), e-mail:hwliuxidian@163.com 
\at   School of Mathematics and Statistics, Xidian University, Xi'an, 710126,   China
}
\date{Received: date / Accepted: date}

\maketitle

\begin{abstract}It is widely accepted  that the stepsize is of great significance to   gradient method. Two efficient gradient methods with approximately optimal stepsizes mainly based on     regularization models  are proposed for   unconstrained optimization. More exactly,  if the objective function is not close to a quadratic function on the line segment between the current and latest iterates,    regularization models are exploited carefully to generate   approximately optimal stepsizes.  Otherwise,   quadratic approximation models are  used. In particular, when the curvature is non-positive,   special   regularization models are developed.  The convergence  of the proposed methods is established under the weak conditions. Extensive numerical experiments indicated the proposed method is superior to the BBQ method (SIAM J. Optim. 2021,31(4), 3068–3096) and other efficient gradient methods,  and is competitive to two  famous and efficient conjugate gradient  software packages CG$ \_ $DESCENT (5.0)  (SIAM J. Optim.   16(1), 170-192, 2005) and CGOPT (1.0) (SIAM J. Optim.  23(1), 296-320, 2013) for   the CUTEr library. Due to the    surprising  efficiency, we believe that gradient methods with approximately optimal stepsizes can become   strong candidates for 	 large-scale unconstrained  optimization.
 
 
 
 
\end{abstract}
\keywords{ Approximately optimal stepsize.  Gradient method.  Global convergence.     Regularization method.  Barzilai-Borwein (BB) method  }\subclass{90C06 \and 65K}
\section{Introduction}

We consider the unconstrained optimization problem:
\begin{align}\label{eq:UnconstrainPro}
\mathop {\min }\limits_{x \in {R^n}} f\left( x \right),
\end{align}
where $f:{R^n} \to R$ is continuously differentiable and its gradient is denoted by $g(x)$.  The gradient method for solving \eqref{eq:UnconstrainPro}  has the form
\begin{align}
{x_{k + 1}} = {x_k} - {\alpha _k}{g_k},
\end{align}
where  ${\alpha _k}$ is the stepsize and ${g_k} = \nabla f({x_k})$. Throughout this paper, ${f_k} = f({x_k})$,  ${s_{k - 1}} = {x_k} - {x_{k - 1}}$, ${y_{k - 1}} = {g_k} - {g_{k - 1}}$   and $\left\| . \right\|$ denotes the Euclidean norm.

It is widely accepted that the stepsize is of great significance to the
theory and   numerical performance of gradient method,  and the stepsize for gradient method has attracted  extensive attentions.
The classical steepest descent method \cite{1}, in which the stepsize is   given by
$
\alpha _k^{SD} = \arg \mathop {\min }\limits_{\alpha  > 0} f({x_k} - \alpha {g_k}),
$
is badly affected by ill conditioning and thus converges slowly \cite{Akaike1959On}. In 1988, Barzilai and Borwein \cite{2} proposed a new   gradient method (BB method), where the famous stepsize (BB stepsize) is given by
\begin{align} \label{eq:BBstepsize}
\alpha _k^{B{B_1}} = \frac{{{{\left\| {{s_{k - 1}}} \right\|}^2}}}{{s_{k - 1}^T{y_{k - 1}}}}\quad\ \text{or}\quad\ \alpha _k^{B{B_2}} = \frac{{s_{k - 1}^T{y_{k - 1}}}}{{{{\left\| {{y_{k - 1}}} \right\|}^2}}}.
\end{align}

\noindent Due to the simplicity and nice numerical efficiency, the BB method has received extensive attentions. The BB method has been shown to be globally \cite{3} and R-linearly \cite{4} convergent for any dimensional strictly convex quadratic functions. In 2021, Li and
Sun \cite{SunRlinearBB2021} presented an interesting improved R-linear convergence result of the BB method.  Dai et al.\cite{5} presented an efficient gradient method by adaptively choosing the BB stepsizes. Raydan\cite{6} proposed a  global BB method by incorporating the nonmonotone line search (GLL line search)\cite{7}. Dai et al.\cite{8}  viewed the BB stepsize from a new angle and  constructed a quadratic model and a conic model   to derive two step sizes for BB-like methods. Based on  a fourth order conic model and some modified secant equations, Biglari and Solimanpur\cite{9} presented some   BB-like methods.  More   BB-like methods can be found in	\cite{ Dai2006CBB,Xiao2010Notes,Nosratipour2017An,Miladinovic2011Scalar}.

 In 2018, Liu et al.\cite{11} viewed the stepsize $\alpha _k^{BB_1}$ from an  approximation model and  introduced a new type of stepsize called approximately optimal stepsize for gradient method.

\noindent \textbf{Definition 1.1}\cite{11} Suppose $f$ is continuously differentiable, and let $\phi_k (\alpha )$ be an approximation model of $f({x_k} - \alpha {g_k}).$ A positive constant ${\alpha_k ^{AOS}}$ is called \textbf{approximately optimal stepsize} associated to $\phi_k (\alpha )$ for gradient method if $\alpha_k ^{AOS}$ satisfies
\begin{align}\label{eq:DefAOS}
{\alpha_k ^{AOS}} = \arg \mathop {\min }\limits_{\alpha  > 0} \phi_k (\alpha ).
\end{align}

Based on \eqref{eq:DefAOS}, it is easy to  obtain    the following simple facts:

 (i)If ${\phi _k}\left( \alpha  \right) =f\left( x_k -\alpha g_k\right) $, then the resulted approximately optimal stepsize corresponds to Cauchy stepsize or optimal stepsize. This is   the reason that we call the stepsize \eqref{eq:DefAOS}   approximately optimal stepsize. 
 
 (ii)If ${\phi _k}\left( \alpha  \right) = {f_k} - \alpha {\left\| {{g_k}} \right\|^2} + \frac{1}{2}{\alpha ^2}g_k^T\left( {\frac{{s_{k - 1}^T{y_{k - 1}}}}{{{{\left\| {{s_{k - 1}}} \right\|}^2}}}I} \right){g_k}$, then the resulted approximately optimal stepsize corresponds to the BB stepsize $\alpha _k^{BB_1}$.  
 
 (iii)If ${\phi _k}\left( \alpha  \right) = {f_k} - \alpha {\left\| {{g_k}} \right\|^2} + \frac{1}{2}{\alpha ^2}g_k^T\left( \dfrac{1}{t} I \right){g_k}$, where $t>0$, then the resulted approximately optimal stepsize corresponds to the fixed stepsize $t$. In fact, for any existing stepsize  $\alpha_k>0$, let  ${\phi _k}\left( \alpha  \right) = {f_k} - \alpha {\left\| {{g_k}} \right\|^2} + \frac{1}{2}{\alpha ^2}g_k^T\left( \dfrac{1}{\alpha_k} I \right){g_k}$, we can easily   see that the resulted approximately optimal stepsize is exactly $\alpha_k$.  
 
 Therefore, all existing stepsizes for gradient methods can be regarded as approximately optimal stepsizes in this sense.  Some   gradient methods with approximately optimal stepsizes \cite{12,Liu2018Several} were   proposed, and the numerical experiments in  \cite{12,Liu2018Several} indicated that these   gradient methods  are very efficient. Gradient methods with approximately optimal stepsizes have illustrated  powerful potentiality for unconstrained optimization.

Besides, an new and important advance for gradient method is the BBQ method \cite{BBQ-HDL}. Motivated by Yuan's  stepsize \cite{Yuan-stepsize}, Huang, Dai and Liu \cite{BBQ-HDL} equiped the Barzilai and Borwein (BB) method with two dimensional quadratic termination property and   proposed a novel stepsize for gradient method (BBQ, corresponding to Algorithm 3.1 in \cite{BBQ-HDL}) for general unconstrained optimization. 

\textbf{Contributions.} According to Definition 1.1, it is not difficult to see that   the effectiveness of approximately optimal stepsize relies heavily on  the approximation model $\phi_k (\alpha )$. To obtain more efficient gradient methods with approximately optimal stepsizes, one should take full advantage of the properties of $f$ at $x_k$ to exploit suitable approximation models including quadratic models and non-quadratic models for deriving approximately optimal stepsize.   Two  efficient gradient methods with approximately optimal stepsizes   are proposed  for   unconstrained optimization in this paper. In the proposed methods, if the objective function $f$ is not close to a quadratic   on the line segment between $x_{k-1}$ and $x_k$,  some  regularization models are exploited to generate   approximately optimal stepsizes. Otherwise, a quadratic approximation model is used to derive approximately optimal stepsize. In particular,   when $s_{k - 1}^T{y_{k - 1}} \le 0,$ some special   regularization models are developed carefully.   The global convergence    of the proposed methods is analyzed. Some numerical results indicate  that the proposed method is  superior to the BBQ method\cite{BBQ-HDL}     and other efficient gradient methods, and is competitive to two  famous conjugate gradient  software packages CGOPT (1.0) \cite{Dai2014A}  and CG$ \_ $DESCENT (5.0) \cite{Hager2005A} for the 145 test problems in the CUTEr library \cite{Gould2003CUTEr}, and  has significant improvement over   CGOPT (1.0) \cite{Dai2014A}  and CG$ \_ $DESCENT (5.0) \cite{Hager2005A} for the  80 test problems mainly from \cite{Andrei2008collection}. It is noted that CGOPT and  CG$ \_ $DESCENT   are  widely treated as two most efficient conjugate gradient software packages.

The rest of the paper is organized as follows. In Section 2, some approximation models including   regularization models and quadratic models are exploited  to generate approximately optimal stepsizes for gradient methods. In Section 3, two efficient gradient methods with the approximately optimal stepsizes are described. The global convergence  of the proposed methods is analyzed in Section 4. In Section 5,  the numerical results are presented. Conclusion  and discussion   are given in the last section.

\section{Derivation of   Approximately Optimal Stepsizes}
Based on the properties of $ f $ at the current iterate $x_k$,    some  approximation models including   regularization  models  and   quadratic models are exploited carefully to derive approximately optimal stepsizes  for   gradient methods in the  section.

As mentioned above,    the effectiveness of approximately optimal stepsize     relies heavily on   approximation model. So  we   take full advantage of the properties of $ f $ at   $x_k$  to construct suitable approximation models for generating approximately optimal stepsizes. The choices of approximation models  are based on the following observations.

Define
\begin{align}
{\mu _k} = \left| {\frac{{2({f_{k - 1}} - {f_k} + g_k^T{s_{k-1}})}}{{s_{k - 1}^T{y_{k - 1}}}} - 1} \right|.
\end{align}
According to \cite{11}, ${\mu _k}$ is an important criterion for judging the degree of $f$ to approximate   quadratic model.   If the condition \cite{8,12}
\begin{align}\label{eq:Quadjus}
{\mu _k} \le {c_1}\quad \text{or} \quad \max \left\{ {{\mu _k},{\mu _{k - 1}}} \right\} \le {c_2},
\end{align}
where    $0 < {c_1} < {c_2},$ holds,   then $f$ might be close to a quadratic function  on the line segment between $x_{k-1}$ and $x_k$.

 When  $f$ is close to a quadratic on the line segment between ${x_{k - 1}}$ and ${x_k},$   quadratic approximation model is preferable. However, if the objective function $f$ possesses high non-linearity,  then  quadratic models might not work very well\cite{14,15}, so  some non-quadratic approximation models should be considered. In recent years,   regularization algorithms for unconstrained optimization,  which are defined as the standard quadratic model plus a regularization term, have become an alternative to trust region and line search schemes\cite{16}. An adaptive regularization algorithm using cubics   (ARC)   was  proposed by Cartis et al.\cite{16}. The trial step  in ARC algorithm is probably computed by minimizing the following   regularization model:
\begin{align}\label{eq:cubicmodel1}
{m_k}({d}) = f({x_k}) + {g^T_k}{d} + \frac{1}{2}d^T{B_k}{d} + \frac{1}{3}{\sigma _k}{\left\| {{d}} \right\|^3 },
\end{align}
where ${B_k}$ is a symmetric   approximation to the Hessian matrix, ${\sigma _k} >0$ is an adaptive positive parameter which can be viewed as the reciprocal of the trust region radius.  And the numerical results in \cite{17} indicated that ARC algorithm is quite efficient. An alternative approach to compute an approximate minimizer of the cubic model has been recently proposed in \cite{18}. In \cite{19}, a nonmonotone cubic overestimation algorithm has been put forward, which follows the one presented in \cite{20}. In \cite{21}, a new algorithm has been designed by combining the   regularization method with line search and nonmonotone techniques.  All of this indicates that when $f$ is not close to a quadratic on the line segment between ${x_{k - 1}}$ and ${x_k}$,   regularization models might serve better than   quadratic models. Based on the above observations, if  $f$ is not close to a quadratic on the line segment between ${x_{k - 1}}$ and ${x_k}$, then we consider the following    regularization models
\begin{align}\label{eq:cubicmodelT}
	{m_k}({d}) = f({x_k}) + {g^T_k}{d} + \frac{1}{2}d^T{B_k}{d} + \frac{1}{p}{\sigma _k\left( p\right) }{\left\| {{d}} \right\|^p },
\end{align}
where $p=3$ or 4, and  ${\sigma _k}\left( p\right) >0$ is a   regularization parameter relative to $p$, otherwise  we construct quadratic models  to generate approximately optimal stepsizes.

We derive the approximately optimal stepsizes for gradient methods in the following four cases.

\textbf{Case I.} $s_{k - 1}^T{y_{k - 1}} > 0$  holds
and the condition \eqref{eq:Quadjus} does not hold.

\textbf{(i)$\textbf{p=3}$ }

In the case, the objective function $f$ might  be not close to a quadratic  on the line segment between ${x_{k - 1}}$ and ${x_k}$,   we thus consider the   regularization model \eqref{eq:cubicmodelT} with  $ d=-\alpha g_k $ and $p=3$:
\begin{align} \label{{eq:cubicmodel2}}
{\phi _{11}}(\alpha ) = f({x_k}) - \alpha g_k^T{g_k} + \frac{1}{2}{\alpha ^2}g_k^T{B_k}{g_k} + \frac{1}{3}{\alpha ^3}{\sigma _k}\left( 3\right){\left\| {{g_k}} \right\|^3}.
\end{align}
 Given that the computational cost and storage, $ B_k $ is generated by  imposing the modified Broyden-Fletcher-Goldfarb-Shanno (BFGS) update formula\cite{22} on a scalar matrix ${D_k}$:
\begin{align}\label{eq:Bk1}
{B_k} = {D_k} - \frac{{{D_k}{s_{k - 1}}s_{k - 1}^T{D_k}}}{{s_{k - 1}^T{D_k}{s_{k - 1}}}} + \frac{{{{\overline y }_{k - 1}}\overline y _{k - 1}^T}}{{s_{k - 1}^T{{\overline y }_{k - 1}}}},
\end{align}
where    ${\overline y _{k - 1}} = {y_{k - 1}} + \frac{{{{\overline r }_k}}}{{{{\left\| {{s_{k - 1}}} \right\|}^2}}}{s_{k - 1}}$ and ${\overline r _k} = 3{({g_k} + {g_{k - 1}})^T}{s_{k - 1}} + 6({f_{k - 1}} - {f_k})$.  Here we take $ D_k $ as ${D_k} = {\xi _0}\frac{{y_{k - 1}^T{y_{k - 1}}}}{{s_{k - 1}^T{y_{k - 1}}}}I$, where ${\xi _0} \ge 1$.
Since there exists ${\mu _1} \in [0,1]$ such that
\begin{align}
{\overline r _k} = 3(s_{k - 1}^T{y_{k - 1}} - s_{k - 1}^T{\nabla ^2}f({x_{k - 1}} + {\mu _1}{s_{k - 1}}){s_{k - 1}}),
\end{align}
to improve the numerical performance  we restrict ${\overline r _k}$ as
\begin{align}\label{eq:trrk}
{\overline r _k} = \min \left\{ {\max \left\{ {{{\overline r }_k}, - \xi_1 s_{k - 1}^T{y_{k - 1}}} \right\},\xi_1 s_{k - 1}^T{y_{k - 1}}} \right\},
\end{align}
 where $0 < {  \xi_1} < 0.1.$

It is not difficult to obtain the following lemma.

\textbf{Lemma 2.1.} \emph{ Suppose that $s_{k - 1}^T{y_{k - 1}} > 0.$ Then $s_{k - 1}^T{\overline y _{k - 1}} > 0$ and $B_k$ is symmetric and  positive definite.}

By imposing  $\frac{{d{\phi _{11}}}}{{d\alpha }} =0$,    we obtain   the equation:  $  - g_k^T{g_k} +\alpha g_k^T{B_k}{g_k} + {\alpha ^2}{\sigma _k}\left( 3\right){\left\| {{g_k}} \right\|^3}=0 .$  Since
\begin{align}
\Delta_{11}  = {(g_k^T{B_k}{g_k})^2} + 4{\sigma _k}\left( 3\right){\left\| {{g_k}} \right\|^5} > 0,
\end{align}
  by solving the above equation we can obtain the approximately optimal stepsize
\begin{align}\label{eq:AOS1}
\bar \alpha _k^{AOS(11)}= \frac{{2{{\left\| {{g_k}} \right\|}^2}}}{{\sqrt{ \Delta_{11}}}  + g_k^T{B_k}{g_k}}.
\end{align}
where $B_k$ is given by \eqref{eq:Bk1} with \eqref{eq:trrk}.

It is observed by numerical experiments that  the bound $\left[ {\alpha _k^{B{B_2}},\alpha _k^{B{B_1}}} \right]$ for $\bar \alpha_k^{AOS(11)}$ is very preferable. Therefore, if  $s_{k - 1}^T{y_{k - 1}} > 0 $ and the condition \eqref{eq:Quadjus} does not hold, then we take the following  truncated approximately optimal stepsize
\begin{align}\label{eq:AOS11}
\alpha _k^{AOS(11)} = \max \left\{ {\min \left\{ {\bar \alpha _k^{AOS(11)},\alpha _k^{B{B_1}}} \right\} \alpha _k^{B{B_2}}} \right\} 
\end{align}
for gradient method.

$\textbf{(ii)p=4}$

We   consider the   regularization model \eqref{eq:cubicmodelT} with  $ d=-\alpha g_k $ and $p=4$:
\begin{align} \label{{eq:cubicmode22}}
	{\phi _{12}}(\alpha ) = f({x_k}) - \alpha g_k^T{g_k} + \frac{1}{2}{\alpha ^2}g_k^T{B_k}{g_k} + \frac{1}{4}{\alpha ^4}{\sigma _k}(4){\left\| {{g_k}} \right\|^4},
\end{align}
where $B_k$ is given  \eqref{eq:Bk1} with \eqref{eq:trrk} for the sake of simplicity. 

By imposing $\frac{{d{\phi _{12}}}}{{d\alpha }} =0$, we get the equation   $  - g_k^T{g_k} +\alpha g_k^T{B_k}{g_k} + {\alpha ^3} {\sigma _k}(4){\left\| {{g_k}} \right\|^4} =0.$  Since
\begin{align}
	\Delta_{12}  = \dfrac{1}{4\sigma_k(4)\left\|g_k \right\|^4 }+\dfrac{\left( g_k^TB_kg_k\right)^3 }{ 27{\sigma _k}(4) \left\|g_k \right\|^2 } > 0,
\end{align}
 the above equation only has a real root and two
 imaginary     roots, and thus the approximately optimal stepsize is  the real root:
 \begin{align}\label{eq:barAOS12}\bar \alpha_k^{AOS(12)} {\rm{ = }}\sqrt[3]{{\frac{1}{{2{\sigma _k}(4){{\left\| {{g_k}} \right\|}^2}}}{\rm{ + }}\sqrt {\Delta_{12} } }} + \sqrt[3]{{\frac{1}{{2{\sigma _k}(4){{\left\| {{g_k}} \right\|}^2}}} - \sqrt {\Delta_{12}} }}. \end{align}

 Similar to the case of $p=3$, we also impose the bound   $\left[ {\alpha _k^{B{B_2}},\alpha _k^{B{B_1}}} \right]$ for $\bar \alpha_k^{AOS(12)}$.  Therefore, if  $s_{k - 1}^T{y_{k - 1}} > 0 $ holds and the condition \eqref{eq:Quadjus} does not hold, then we  take the following truncated  approximately  optimal stepsize
 \begin{align}\label{eq:AOS12}
 \alpha _k^{AOS(12)} = \max \left\{ {\min \left\{ {\bar \alpha _k^{AOS(12)},\alpha _k^{B{B_1}}} \right\},\alpha _k^{B{B_2}}} \right\}
 \end{align}
  for the gradient method.

  \indent \textbf{The choice of regularization parameter in the regularization model}

  When the regularization models are applied, the regularization parameter $\sigma_k(p)$ should be determined properly. The regularization parameter is significant to the effectiveness of regularization models.  However, it is  universally acknowledged that it is challenging to determine a proper regularization parameter $\sigma_k\left( p\right) $.  Some ways \cite{Cartis2011Adaptive,Gould2012Updating} were developed to determine the regularization parameter $\sigma_k\left( p\right) $, including the interpolation condition and the trust-region strategy. Here we   use the interpolation condition to determine the regularization parameter $\sigma_k\left( p\right) $:
  \[{f_{k - 1}} = {f_k} - g_k^T{s_{k - 1}} + \frac{1}{2}s_{k - 1}^T{B_k}{s_{k - 1}} + \frac{{{\sigma _k\left(p \right) }}}{p}{\left\| {{s_{k - 1}}} \right\|^{p }},\]
  which implies that
  \[{\sigma _k\left( p\right) } = \frac{{p\left( {{f_{k - 1}} - {f_k} + g_k^T{s_{k - 1}} - \frac{1}{2}\left( {s_{k - 1}^T{y_{k - 1}} + {{\overline r }_k}} \right)} \right)}}{{{{\left\| {{s_{k - 1}}} \right\|}^p}}},\]
 where ${{\overline r }_k} $ is given by \eqref{eq:trrk} and $p=3 \;\text{or} \;4$. To improve the numerical performance and make it to be positive, we take the following truncated form:
  \begin{align}\label{eq:sigma1}{\sigma _k\left( p\right)} = \max \left\{ {\min \left\{ {\left| {\sigma _k\left( p\right)}\right| ,{\sigma _{\max }}} \right\},{\sigma _{\min }}} \right\},  \end{align}
  where $0< \sigma_{\min} < \sigma_{\max}   $ and $p=3 \;\text{or} \;4$.

\textbf{Case II.} $s_{k - 1}^T{y_{k - 1}} > 0$ holds
and the condition \eqref{eq:Quadjus}   holds.

	In the case,  the objective function $f$ might  be   close a quadratic  on the line segment between ${x_{k - 1}}$ and ${x_k}$,   we thus   consider the following quadratic approximation model:
	\begin{align} \label{eq:Quadmodel4}
		{\phi _2}(\alpha ) = f({x_k}) - \alpha g_k^T{g_k} + \frac{1}{2}{\alpha ^2}g_k^T{B_k}{g_k},
	\end{align}
	where ${B_k}$ is    given by  \eqref{eq:Bk1} with \eqref{eq:trrk}  for simplicity. {\tiny }
	By imposing $ \dfrac{d\phi_2}{d \alpha}=0 $, we can easily obtain the approximately optimal stepsize:
	\begin{align}\label{eq:AOS4}
		 \alpha _k^{AOS(2)} =\dfrac{g_k^T{g_k}}{g_k^TB_k{g_k}}.
	\end{align}
 	It is also observed by numerical experiments that the bound $[\alpha _k^{B{B_2}},\alpha _k^{B{B_1}}]$ for $\bar \alpha _k^{AOS(2)}$ is very preferable. Therefore, if $s_{k - 1}^T{y_{k - 1}} > 0$ holds and the condition \eqref{eq:Quadjus} holds, then we take the truncated approximately optimal stepsize
	\begin{align}\label{eq:AOS33}
		\alpha _k^{AOS(2)} = \max \left\{ {\min \left\{ {  \alpha _k^{AOS(2)},\alpha _k^{B{B_1}}} \right\},\alpha _k^{B{B_2}}} \right\}
	\end{align}
for gradient method.

\textbf{Case III.}  $s_{k - 1}^T{y_{k - 1}}  \le 0$ holds and the condition \eqref{eq:conditionfor2}  holds 	

When $s_{k - 1}^T{y_{k - 1}} \le 0$, $ f $ may   enjoy poor  properties at some neighbors of $ x_k $, it is  thus difficult to determine suitable stepsize for gradient method. In some modified BB methods\cite{8,9},  the initial stepsize is usually   set simply  to ${\alpha _k} = 10^{30}$ when   $s_{k - 1}^T{y_{k - 1}} \le 0$. As a result, it will cause large computational cost for seeking a suitable stepsize in a line search for gradient method.


%
%

It follows from    $s_{k - 1}^T{y_{k - 1}} \le 0$ that $
0<\frac{{\left\| {{g_{k - 1}}} \right\|}}{{\left\| {{g_k}} \right\|}} \le 1.$   Consequently, if the following condition:
   \begin{align} \label{eq:conditionfor2}\frac{{{{\left\| {{g_{k - 1}}} \right\|}^2}}}{{{{\left\| {{g_k}} \right\|}^2}}} \ge {\xi _2}, \end{align}
   where $0< {\xi _2} <1$   is close  to 1, holds,  then ${g_k}$ and ${g_{k - 1}}$   incline to be collinear and are approximately equal. Based on the above observation,  we will give a new way to estimate   $  g_k^T{\nabla ^2}f({x_k}){g_k}   $ in approximation model. Therefore, when $s_{k - 1}^T{y_{k - 1}} \le 0$, we    construct a regularization   model to derive approximately optimal stepsizes based on the condition \eqref{eq:conditionfor2}.


 $\textbf{(i)p=3}$

  Suppose for the moment that $ f $ is twice continuously differentiable, we consider the following regularization
  model:
  \begin{align} \label{eq:cubicmodel20}
  \phi (\alpha ) = {f_k} - \alpha g_k^T{g_k} + \frac{1}{2}{\alpha ^2}g_k^T{\nabla ^2}f({x_k}){g_k} + \frac{{{\sigma _k}(3)}}{3}{\alpha ^3}{\left\| {{g_k}} \right\|^3}.
  \end{align}
 When the condition \eqref{eq:conditionfor2}  holds, we use  $ g_{k - 1}^T{\nabla ^2}f({x_k}){g_{k - 1}}  $ to approximate $ g_k^T{\nabla ^2}f({x_k}){g_k}  $ and thus get that
\begin{align}\label{eq:gBg}
g_k^T{\nabla ^2}f({x_k}){g_k} \approx g_{k - 1}^T{\nabla ^2}f({x_k}){g_{k - 1}} \approx \frac{{\left| {{{\left( {g({x_k} + {\alpha _{k - 1}}{g_{k - 1}}) - g({x_k})} \right)}^T}{g_{k - 1}}} \right|}}{{{\alpha _{k - 1}}}} = \frac{{\left| {s_{k - 1}^T{y_{k - 1}}} \right|}}{{\alpha _{k - 1}^2}},
\end{align}
which gives    the following approximation model:
\begin{align*}
{ {\phi}  _{31}}(\alpha ) = {f_k} - {\alpha}g_k^T{g_k} + \frac{1}{2}{\alpha ^2}\frac{{\left| {s_{k - 1}^T{y_{k - 1}}} \right|}}{{\alpha _{k - 1}^2}} + \frac{{{\sigma _k(3)}}}{3}{\alpha ^3}{\left\| {{g_k}} \right\|^3}.
\end{align*}
By imposing $\frac{{d{ {\phi}  _{31}}}}{{d\alpha }} = 0$, we get the equation $ - {\left\| {{g_k}} \right\|^2} + \alpha \frac{{\left| {s_{k - 1}^T{y_{k - 1}}} \right|}}{{\alpha _{k - 1}^2}} + {\alpha ^2}{\sigma _k(3)}{\left\| {{g_k}} \right\|^3} = 0$. Since $${\Delta _{31}} = \frac{{{{\left| {s_{k - 1}^T{y_{k - 1}}} \right|}^2}}}{{\alpha _{k - 1}^4}} + 4{\sigma _k(3)}{\left\| {{g_k}} \right\|^5} > 0,$$      the above equation only has a real root and two
imaginary     roots, and thus the approximately optimal stepsize is  the real root:
\begin{align}\label{eq:stepsize31}
\alpha _k  ^{AOS(31)}= \frac{{2{{\left\| {{g_k}} \right\|}^2}\alpha _{k - 1}^2}}{{\sqrt {{{\left| {s_{k - 1}^T{y_{k - 1}}} \right|}^2} + 4\alpha _{k - 1}^4{\sigma _k(3)}{{\left\| {{g_k}} \right\|}^5}}  + \left| {s_{k - 1}^T{y_{k - 1}}} \right|}}.
\end{align}

$\textbf{(ii)p=4}$

Suppose for the moment that $ f $ is twice continuously differentiable, we consider the following regularization  model:
\begin{align} \label{eq:cubicmodel20}
	\phi (\alpha ) = {f_k} - \alpha g_k^T{g_k} + \frac{1}{2}{\alpha ^2}g_k^T{\nabla ^2}f({x_k}){g_k} + \frac{{{\sigma _k(4)}}}{4}{\alpha ^4}{\left\| {{g_k}} \right\|^4}.
\end{align}

Using \eqref{eq:gBg}, we   get the following model:
\begin{align*}
	{ {\phi}  _{32}}(\alpha ) = {f_k} - {\alpha}g_k^T{g_k} + \frac{1}{2}{\alpha ^2}\frac{{\left| {s_{k - 1}^T{y_{k - 1}}} \right|}}{{\alpha _{k - 1}^2}} + \frac{{{\sigma _k(4)}}}{4}{\alpha ^4}{\left\| {{g_k}} \right\|^4}.
\end{align*}

\noindent By imposing  $\frac{{d{\phi _{32}}}}{{d\alpha }} =0$, we obtain the equation: $ {\sigma _k(4)}{\left\| {{g_k}} \right\|^4}{\alpha ^3} + \frac{{\left| {s_{k - 1}^T{y_{k - 1}}} \right|}}{{\alpha _{k - 1}^2}}{\alpha ^2} - {\left\| {{g_k}} \right\|^2} = 0.$
Since
\[\Delta_{32}  = \frac{1}{{4\sigma _k^2(4){{\left\| {{g_k}} \right\|}^4}}} + \frac{{{{\left| {s_{k - 1}^T{y_{k - 1}}} \right|}^3}}}{{27\sigma _k^3(4)\alpha _{k - 1}^6{{\left\| {{g_k}} \right\|}^{12}}}}>0,\]
 the above equation only has a real root and two
 imaginary     roots, and thus the approximately optimal stepsize is  the real root:
\begin{align}\label{eq:AOS32}\alpha _k^{AOS(32)} = \sqrt[3]{{\frac{1}{{2{\sigma _k(4)}{{\left\| {{g_k}} \right\|}^2}}} + \sqrt {\Delta_{32}} }} + \sqrt[3]{{\frac{1}{{2{\sigma _k(4)}{{\left\| {{g_k}} \right\|}^2}}} - \sqrt {\Delta_{32}} }}.\end{align}

  \textbf{The choice of regularization parameter in the regularization model}

Similar to Case I, we also use the interpolation condition to determine the regularization parameter $\sigma_k\left( p\right) $:
\[{f_{k - 1}} = {f_k} - g_k^T{s_{k - 1}} + \frac{1}{2}s_{k - 1}^T{y_{k - 1}} + \frac{{{\sigma _k\left(p \right) }}}{p}{\left\| {{s_{k - 1}}} \right\|^{p }},\]
which implies that
\[{\sigma _k\left( p\right) } = \frac{{p\left( {{f_{k - 1}} - {f_k} + g_k^T{s_{k - 1}} - \frac{1}{2} {s_{k - 1}^T{y_{k - 1}}  }  } \right)}}{{{{\left\| {{s_{k - 1}}} \right\|}^{p }}}}.\]
Here $p=3 \; \text{or} \; 4$.
To improve the numerical performance and make it to  be positive, we take the following truncation form:
\begin{align}\label{eq:sigma2}{\sigma _k\left( p\right)} = \max \left\{ {\min \left\{ {\left| {\sigma _k\left( p\right)}\right| ,{\sigma _{\max }}} \right\},{\sigma _{\min }}} \right\}, \end{align}
where $0< \sigma_{\min} < \sigma_{\max}   $ are the same as that in \eqref{eq:sigma1} and $p=3 \; \text{or} \; 4$.

\textbf{Case IV.}  $s_{k - 1}^T{y_{k - 1}}  \le 0$ holds and the condition \eqref{eq:conditionfor2}  does not hold

It also has been shown that if $\alpha_k^{\text{BB}}$   is reused in a cyclic fashion, then the convergence rate is accelerated \cite{Miladinovi2011}. It appears that the stepsize $\alpha_{k-1}$ may
provide some important information for the current stepsize. As a result,
we  take $\alpha_k= \xi_3 \alpha_{k-1}$ as the stepsize,   where $\xi_3 >0$. In actual, the stepsize can also be regarded as the approximately optimal stepsize. By taking $B_k =\frac{1}{\xi_3\alpha_{k-1}}I  $,   we can get the following quadratic approximation model
	\begin{align} \label{eq:Quadmodel5}
	{\phi _4}(\alpha ) = f({x_k}) - \alpha g_k^T{g_k} + \frac{1}{2}{\alpha ^2}g_k^T\left(\dfrac{1}{\xi_3\alpha_{k-1}}I \right) {g_k}.
	\end{align}
	By imposing $\frac{{d{ {\phi}  _4}}}{{d\alpha }} = 0$, we obtain the approximately optimal stepsize:
	\begin{align}\label{eq:stepsize4}
	\alpha _k  ^{AOS(4)}=\xi_3 \alpha_{k-1}.
	\end{align}

\section{Two Efficient Gradient Methods with Approximately Optimal Stepsizes}

We   describe two efficient  gradient methods with approximately optimal stepsizes in the section.

The famous nonmonotone line search (GLL line search) \cite{7} was firstly  incorporated into the BB method \cite{6}. Though  GLL line search works well in many cases, there are some drawbacks, for example, some good
function values may be discarded, or the numerical performance depends very much
on the choice of a pre-fixed memory constant.
 To overcome the above drawbacks, another well-known nonmonotone Armijo line search (Zhang-Hager line search)  \cite{10} was proposed by Zhang and Hager and  is  defined as
		\begin{align} \label{eq:ZhangHageLS1}
			f({x_k} - {\alpha  }{g_k}) \le {C_k} - \delta {\alpha  }{\left\| {{g_k}} \right\|^2},
		\end{align}
where $ 0<\delta<1 $,		\begin{align}\label{eq:ZhangHageLS2}
Q_0=1,\;\;	{Q_{k + 1}} = {\eta _k}{Q_k} + 1,\quad C_0=f(x_0), \;\;  {C_{k + 1}} = {{({\eta _k}{Q_k}{C_k} + f({x_{k + 1}}))} \mathord{\left/
			{\vphantom {{({\eta _k}{Q_k}{C_k} + f({x_{k + 1}}))} {{Q_{k + 1}}}}} \right.
			\kern-\nulldelimiterspace} {{Q_{k + 1}}},\;\;0<\eta_k\le 1}.
\end{align}
 It is observed that     Zhang-Hager line search \cite{10}  is usually preferable for the BB-like methods. To improve the numerical performance and obtain nice convergence, we take $ \eta_k $ as :
		\begin{align}\label{eq:ZhangHageLS3}
			{\eta _k} = \left\{ \begin{array}{l}
				c,\quad \bmod (k,n) = n - 1,\\
				1,\quad \bmod (k,n) \ne n - 1,
			\end{array} \right.			
		\end{align}	
	where $0 < c < 1$ and $\bmod (k,n)$ represents the residue for $k$ modulo $n$. As a result,   Zhang-Hager line search \cite{10}  with \eqref{eq:ZhangHageLS3} and the following strategy
 \cite{23}:
\begin{align}\label{eq:ZhangHageLS4}
{\alpha  } = \left\{ \begin{array}{l}
{\overline \alpha   },\quad \quad {\rm{if}}\;{\alpha } > 0.1\alpha^{(0)} _k\;{\rm{and}}\;{\overline \alpha   } \in [0.1\alpha^{(0)} _k,0.9{\alpha  }],\\
0.5{\alpha  },\;\;\;{\rm{otherwise}},
\end{array} \right.
\end{align} 	
	where $\alpha^{(0)}$ is   approximately optimal stepsize described in Section 2 and  ${\overline \alpha  }$ is obtained by a quadratic interpolation at $ x_k $ and ${x_k} - {\alpha  }{g_k},$	  is used in the proposed methods.\\

We   describe the  gradient method with approximately optimal stepsize   (GM\_AOS (Reg p=3)) in detail.

\begin{algorithm}\label{Algorithm1}
	\centering
	\caption{ GM\_AOS (Reg p=3)}
	\begin{algorithmic}
		\STATE\textbf {Step 0.} Initialization. Given    ${x_0} \in {R^n},$  $\varepsilon  > 0,$  $\delta,$ $c,\;c_1,\;c_2$  ${\alpha_{\max }},$ ${\alpha_{\min }},$ $\alpha^0 _0,$
	 ${\sigma _{\min }},$   ${\sigma _{\max }},$  ${\xi _0}, \;{\xi _1}, \;{\xi _2},\; \xi_3$.     Set
	
\quad \quad ${Q_0} = 1,$ ${C_0} = {f_0}$  and $k = 0.$\\
		\STATE\textbf {Step 1.} If $\left\| {{g_k}} \right\|_{\infty} \le \varepsilon ,$ then stop. \\
		\STATE\textbf {Step 2.} Compute   approximately optimal  stepsize.
		
	\quad	\textbf{  2.1} If $k = 0$, then set  $ \alpha=\alpha_0^{(0)} $
	and go to Step 3.
	
	\quad	\textbf{  2.2}  If $s_{k - 1}^T{y_{k - 1}} > 0$ holds and  the condition \eqref{eq:Quadjus}   does not  hold,    then compute ${\alpha _k}$ by \eqref{eq:AOS11} and update
	
\quad \quad \quad	$\sigma_k$ by $\eqref{eq:sigma1}$ with $p=3$.
	
	\quad	\textbf{  2.3}  If   $s_{k - 1}^T{y_{k - 1}} > 0$ holds and the condition \eqref{eq:Quadjus}  holds,    then compute ${\alpha _k}$ by \eqref{eq:AOS33}.

	\quad	\textbf{  2.4}  If $s_{k - 1}^T{y_{k - 1}} \le 0$ holds and     the condition \eqref{eq:conditionfor2}    holds,    then compute ${\alpha _k}$ by \eqref{eq:stepsize31} and update
	
	\quad \quad \quad $\sigma_k$ by $\eqref{eq:sigma2}$ with $p=3$.
	
	\quad	\textbf{  2.5}  If $s_{k - 1}^T{y_{k - 1}} \le 0$ holds and     the condition \eqref{eq:conditionfor2} does not hold,    then compute ${\alpha _k}$ by \eqref{eq:stepsize4}.

	\quad	\textbf{  2.6} Set	 ${\alpha^{(0)} _k} = \max \left\{ {\min \left\{ {{\alpha _k},{\alpha _{\max }}} \right\},{\alpha _{\min }}} \right\}$ and $\alpha= {\alpha^{(0)} _k}$.
		\STATE\textbf {Step 3}. Line search. If \eqref{eq:ZhangHageLS1} holds,	then go to Step 4, otherwise  update  $ \alpha  $ by \eqref{eq:ZhangHageLS4}  and  go to		 Step 3.  \\
		\STATE\textbf {Step 4.}  Update ${Q_{k + 1}},\;{C_{k + 1}} $ and $ \eta_k $ by \eqref{eq:ZhangHageLS2} and  \eqref{eq:ZhangHageLS3}.
		\STATE\textbf {Step 5.}  Set $\alpha_k =\alpha, \; {x_{k + 1}} = {x_k} - {\alpha _k}{g_k},$ $k = k + 1,$ and go to Step 1.
	\end{algorithmic}
\end{algorithm}

\noindent\textbf{Remark.} If ``\textbf{2.2} If $s_{k - 1}^T{y_{k - 1}} > 0$ holds  and  the condition \eqref{eq:Quadjus}   does not hold,    then compute ${\alpha _k}$ by \eqref{eq:AOS12} and update  	$\sigma_k$ by $\eqref{eq:sigma1}$ with $p=4$" and ``\textbf{2.4}  If $s_{k - 1}^T{y_{k - 1}} \le 0$ holds and     the condition \eqref{eq:conditionfor2}  holds,    then compute ${\alpha _k}$ by \eqref{eq:AOS32} and update  $\sigma_k$ by $\eqref{eq:sigma2}$ with $p=4$"     are used to replace   of   \textbf{2.2} and \textbf{2.4}   of  Algorithm 1, relatively, then the resulting method corresponds to  another gradient method with approximately optimal stepsize called GM$ \_ $AOS (Reg  p=4). We use GM$ \_ $AOS (Reg)
to denote either GM$ \_ $AOS (Reg p=3) or GM$ \_ $AOS (Reg  p=4).

\section{Convergence Analysis}

In the section    the  global convergence  of  GM$\_$AOS (Reg)  is analyzed under weak  conditions. In the  convergence analysis   the following assumptions are done. \\
D1. $f(x)$ is continuously differentiable on $\mathbb R^n$. \\
D2. $f(x)$ is bounded below on $\mathbb R^n$. \\
D3. The gradient $g(x)$    is \textbf{ uniformly    continuous} on  $ \mathbb   R^n $.

\begin{lemma} 	\label{thm lemma5-1}  For $ Q_k $ in \eqref{eq:ZhangHageLS2} , we have ${Q_{k + 1}} \le 1 + \dfrac{n}{{1 - c}}$.
\end{lemma}

\noindent	 \textbf{ Proof}
It follows from  (\ref{eq:ZhangHageLS2})   that \[{Q_{k + 1}} = 1 + \sum\limits_{j = 0}^k {\prod\limits_{i = 0}^j {{\eta _{k - i}}} },   \]
which together with (\ref{eq:ZhangHageLS3}) suggests that \begin{equation}\label{eq:ConverAnNewQk}{Q_{k + 1}} = \left\{ {\begin{array}{*{20}{c}}
	{1 + n\sum\limits_{i = 1}^{\left( {k + 1} \right)/n} {{c^i}} \;,\;\;\;\;\;\;\;\;\;\;\;\;\;\;\;\;\;\;\;\;\;\;\;\;\;\;\;\;\;\;\;\;\;\;\;\;\text{if}\;\bmod (k,n) = n - 1,} \\
	{1 + \left( {1 + \bmod \left( {k,n} \right) + n\sum\limits_{i = 1}^{\left\lfloor {k/n} \right\rfloor } {{c^i}} } \right)\;,\;\;\;\;\;\text{if}\;\bmod (k,n) \ne n - 1,}
	\end{array}} \right. \end{equation}
where $ \left \lfloor \cdot  \right\rfloor  $ is the floor function.

By (\ref{eq:ConverAnNewQk}) and the fact that $0< c<1 $, we obtain  that
\[{Q_{k + 1}} \leqslant 1 + \left( {n + n\sum\limits_{i = 1}^{\left\lfloor {k/n} \right\rfloor  + 1} {{c^i}} } \right) \leqslant 1 + \left( {n + n\sum\limits_{i = 1}^{k + 1} {{c^i}} } \right) = 1 + n\sum\limits_{i = 0}^{k + 1} {{c^i}}  = 1 + \frac{{n\left( {1 - {c^{k + 2}}} \right)}}{{1 - c}} \leqslant 1 + \frac{n}{{1 - c}},  \]
which     completes the proof.                  \qed	


\begin{lemma} \label{thm lemma5-2}
	Suppose that D1, D2 and D3 hold.  Then, \begin{equation}\label{eq:ConverAnwkfkCkCk1}
	{f_{k + 1}}  \le  {C_{k + 1}}  \le   {C_k}.  \end{equation}
\end{lemma}
\textbf{Proof}  According to (\ref{eq:ZhangHageLS1}) and  (\ref{eq:ZhangHageLS2}), we have
\[{C_{k + 1}} = \frac{{{\eta _k}{Q_k}{C_k} + f_{k + 1}}}{{{Q_{k + 1}}}} = {C_k} + \frac{{f_{k + 1} - {C_k}}}{{{Q_{k + 1}}}} \le {C_k} \]
and
\[{C_{k + 1}} = \frac{{{\eta _k}{Q_k}{C_k} + {f_{k + 1}}}}{{{Q_{k + 1}}}} = \frac{{{\eta _k}{Q_k}}}{{{\eta _k}{Q_k} + 1}}{C_k} + \frac{1}{{{\eta _k}{Q_k} + 1}}{f_{k + 1}} \ge \frac{{{\eta _k}{Q_k}}}{{{\eta _k}{Q_k} + 1}}{f_{k + 1}} + \frac{1}{{{\eta _k}{Q_k} + 1}}{f_{k + 1}} = {f_{k + 1}}.\]
As a result, the inequality (\ref{eq:ConverAnwkfkCkCk1}) holds. The proof is completed.                     \qed

The above lemma implies that the sequence $ \left\lbrace C_k \right\rbrace  $ is convergent.

\begin{theorem} \label{thm5-3} Suppose that D1, D2 and D3 hold, and let $ \left\lbrace x_k \right\rbrace  $ be the sequence generated by  GM$\_$AOS (Reg).  Then,
	\begin{equation}\label{eq:ConverAnaResult}\mathop {\lim }\limits_{k \to \infty } \left\|  g_k \right\|  = 0. \end{equation}	
\end{theorem}	
\textbf{ Proof}   By (\ref{eq:ZhangHageLS1}) and (\ref{eq:ZhangHageLS2}), we obtain that
\[{C_{k + 1}} = {C_k} + \frac{{{f_{k + 1}} - {C_k}}}{{{Q_{k + 1}}}} \le {C_k} - \frac{{\sigma {\alpha _k}{{\left\| {{g_k}} \right\|}^2}}}{{{Q_{k + 1}}}},\]
which together with      Lemma \ref{thm lemma5-1} implies that
\begin{equation}\label{eq:ConverAnakgkCk}\frac{\sigma }{{1 + n/(1 - c)}}{\alpha _k}{\left\| {{g_k}} \right\|^2} \le \frac{{\sigma {\alpha _k}{{\left\| {{g_k}} \right\|}^2}}}{{{Q_{k + 1}}}} \le {C_k} - {C_{k + 1}}. \end{equation}
It then follows from Lemma \ref{thm lemma5-2}  and D2 that
\begin{equation}\label{eq:ConverAnwkfkakgk} \mathop {\lim }\limits_{k \to \infty } {\alpha _k}{\left\| {{g_k}} \right\|^2} = 0.  \end{equation}

We  suppose, by way of contradiction, that  there exists a subsequence $ \left\lbrace x_{k_j} \right\rbrace  $ such that
\begin{equation}\label{eq:ConverAnatradition}  \mathop {\lim \;}\limits_{j \to \infty } \left\| {  {g_{{k_j}}}} \right\| = l > 0.  \end{equation} We denote $$\overline \varepsilon   = \left\{ {\begin{array}{*{20}{c}}
	{l/2,\;\;\;\;\;\;{\rm{if}}\;l <  + \infty ,}\\
	{1/2,\;\;\;\;\;{\rm{otherwise}}{\rm{.}}}
	\end{array}} \right.$$
It follows from (\ref{eq:ConverAnatradition}) that there exists a positive integer  $ {j_0} $ such that   \begin{equation}\label{eq:ConverAnagkgk0}  \left\| {  {g_{{k_j}}}} \right\| > \overline \varepsilon, \;\;\;\forall  j>j_0. \end{equation}   
Therefore, we obtain from (\ref{eq:ConverAnwkfkakgk}) that $ \mathop {\lim }\limits_{j \to    \infty } {\alpha _{{k_j}}} = 0$  and
\begin{equation}\label{eq:ConverAnaak0akgk0}  \;\mathop {\lim }\limits_{j \to    \infty } \alpha _{{k_j}}^2{\left\| {{g_{{k_j}}}} \right\|^2} = 0. \end{equation}

By (\ref{eq:ZhangHageLS4}),  we know that there exists  $  \bar \delta _{k_j}  \in [0.1,0.9] $ such that
\begin{equation}\label{eq:ConverAnfkCk}\;f\left( {{x_{{k_j}}} - \frac{{{\alpha _{{k_j}}}}}{ \bar \delta _{k_j} }{g_{{k_j}}}} \right) > {C_{{k_j}}} - \sigma \frac{{{\alpha _{{k_j}}}}}{\bar \delta _{k_j} }{\left\| {{g_{{k_j}}}} \right\|^2}.  \end{equation}
Combining (\ref{eq:ConverAnfkCk}) and   $ {f\left( {{x_{{k_j}}} - {\alpha _{{k_j}}}{g_{{k_j}}}} \right)} \le {C_{{k_j}}} - \sigma {\alpha _{{k_j}}}{\left\| {{g_{{k_j}}}} \right\|^2}$,	we obtain that \[f\left( {{x_{{k_j}}} - \frac{{{\alpha _{{k_j}}}}}{\bar \delta _{k_j} }{g_{{k_j}}}} \right) - f\left( {{x_{{k_j}}} - {\alpha _{{k_j}}}{g_{{k_j}}}} \right) >  - \sigma \left( {\frac{1}{\bar \delta _{k_j} } - 1} \right){\alpha _{{k_j}}}{\left\| {{g_{{k_j}}}} \right\|^2}.\]  	

\noindent{It follows from the mean-value theorem that   there exists  $ {w_{{k_j}}} \in \left[ {0,1} \right] $ such that   	
	\[f\left( {{x_{{k_j}}} - \frac{{{\alpha _{{k_j}}}}}{\bar \delta _{k_j} }{g_{{k_j}}}} \right) - f\left( {{x_{{k_j}}} - {\alpha _{{k_j}}}{g_{{k_j}}}} \right) =  - \left( {\frac{1}{\bar \delta _{k_j} } - 1} \right){\alpha _{{k_j}}}g{\left( {{u_{{k_j}}}} \right)^T}g_{{k_j}},\]
	where ${u_{{k_j}}} = {x_{{k_j}}} - \left[ {1 + {w_{{k_j}}}\left( {1/{{\bar \delta }_{{k_j}}} - 1} \right)} \right]{\alpha _{{k_j}}}{g_{{k_j}}}$.   	
	Therefore, we get  that \[ - \left( {\frac{1}{{\bar \delta _{k_j} }} - 1} \right){\alpha _{{k_j}}}g{\left( {{u_{{k_j}}}} \right)^T}{g_{{k_j}}} >  - \sigma \left( {\frac{1}{{\bar \delta _{k_j} }} - 1} \right){\alpha _{{k_j}}}{\left\| {{g_{{k_j}}}} \right\|^2},\]
	which implies that
	$\;{\left( {{g_{{k_j}}} - g\left( {{u_{{k_j}}}} \right)} \right)^T}\frac{{{g_{{k_j}}}}}{{\left\| {{g_{{k_j}}}} \right\|}} > \left( {1 - \sigma } \right)\left\| {{g_{{k_j}}}} \right\|.$  	
	According to (\ref{eq:ConverAnagkgk0}), we know that   	
	\begin{equation}\label{eq:ConverAngkgk}\left\| {{g_{{k_j}}} - g({u_{{k_j}}})} \right\| \ge \left( {{g_{{k_j}}} - g({u_{{k_j}}})} \right)^T\frac{{{g_{{k_j}}}}}{{\left\| {{g_{{k_j}}}} \right\|}} > \left( {1 - \sigma } \right)\left\| {{g_{{k_j}}}} \right\|\; > \left( {1 - \sigma } \right)\bar \varepsilon, \;\;\; \forall j >  {j_0}. \end{equation}	 	
	According to   (\ref{eq:ConverAnwkfkakgk}), (\ref{eq:ConverAnaak0akgk0}) and $1 \le  {1 + {w_{{k_j}}}\left( {1/{{\bar \delta }_{{k_j}}} - 1} \right)}  \le 10$, we know that
	\begin{equation}\label{eq:ConverAnwkakgk0}\mathop {\lim }\limits_{j \to  + \infty } \left[ {{w_{{k_j}}}\left( {1/{{\bar \delta }_{{k_j}}} - 1} \right) + 1} \right]{\alpha _{{k_j}}}\left\| {{g_{{k_j}}}} \right\| \to 0.   \end{equation}   Since the gradient $ g $ is uniformly continuous,   for $ \frac{{\left( {1 - \sigma } \right) \bar \varepsilon }}{2} $, one can find    $ \zeta >0$ depending only on $ \frac{{\left( {1 - \sigma } \right)\bar \varepsilon }}{2} $ such that
	$ \left\| {{g_{{k_j}}} - g\left( {{u_{{k_j}}}} \right)} \right\| \le \frac{{\left( {1 - \sigma } \right)}}{2}\bar \varepsilon  $ holds whenever $\left\| {{x_{{k_j}}} - {u_{{k_j}}}} \right\| = \left[ {{w_{{k_j}}}\left( {1/{{\bar \delta }_{{k_j}}} - 1} \right) +1} \right]{\alpha _{{k_j}}}\left\| {{g_{{k_j}}}} \right\| < \zeta .$  By  (\ref{eq:ConverAnwkakgk0}),  we know    that there exists an integer $ j_1 >0 $ such that  $$\left\| {{x_{{k_j}}} - {u_{{k_j}}}} \right\| = \left[ {{w_{{k_j}}}\left( {1/{{\bar \delta }_{{k_j}}} - 1} \right) + 1} \right]{\alpha _{{k_j}}}\left\| {{g_{{k_j}}}} \right\| < \zeta   $$ holds  for any $ j>j_1 $. As a result,
	$ \left\| {{g_{{k_j}}} - g\left( {{u_{{k_j}}}} \right)} \right\| \le \frac{{\left( {1 - \sigma } \right)}}{2}\bar \varepsilon   $ holds for any $ j>j_1 $,   which  contradicts  (\ref{eq:ConverAngkgk}) when $j\ge\max\left\lbrace j_0,j_1\right\rbrace $.    Therefore, there no exists a subsequence $ \left\lbrace x_{k_j} \right\rbrace  $ satisfying (\ref{eq:ConverAnatradition}), which implies  (\ref{eq:ConverAnaResult}).
	The proof is completed.  \qed
	

\section{Numerical Experiments}

 We compare  GM$ \_ $AOS (Reg) with  GM\_AOS (1.2) \cite{Liu2018Several}, the BB method,
CGOPT (1.0) \cite{Dai2014A}, CG\_DESCENT (5.0) \cite{Hager2005A}   and BBQ method \cite{BBQ-HDL} (corresponding to Algorithm 3.1 in \cite{BBQ-HDL}) in the section. It is widely accepted  that CGOPT \cite{Dai2014A} and CG\_DESCENT \cite{Hager2005A} are  the two most famous and efficient conjugate gradient software packages.   The codes of the BB method, GM$ \_ $AOS (1.2) \cite{Liu2018Several} and  GM$ \_ $AOS (Reg) were implemented by C language, and the C codes of CG$ \_ $DESCENT (5.0) and CGOPT (1.0) can be downloaded  from     \url{http://users.clas.ufl.edu/hager/papers/Software} and \url{http://coa.amss.ac.cn/wordpress/?page_id=21}, respectively. The matlab code of BBQ can be found in Dai's homepage: \url{http://lsec.cc.ac.cn/~dyh/software.html}.   The C code of  GM$ \_ $AOS (Reg) and the detailed numerical result  will be available in our website finally.  Two test sets were used, which include the 145 test problems in the CUTEr library \cite{Gould2003CUTEr} (we call it   CUTEr145  for short) and the 80 test problems mainly from \cite{Andrei2008collection} collected by Andrei (we call it  Andr80  for short),    respectively. The two test sets can be    found    in Hager's website \url{http://users.clas.ufl.edu/hager/papers/CG/results6.0.txt} and  Andrei's homepage \url{http://camo.ici.ro/neculai/AHYBRIDM}, respectively.  The dimensions of the test problem  in the test set  CUTEr145   are default and the dimension of  each test problems  in the  test set  Andr80 is set to 10,000.
 All numerical experiments were done in Ubuntu 10.04 LTS     in  a VMware Workstation 10.0 installed in Win 10.

 We choose the following  parameters for  GM$ \_ $AOS (Reg): $\varepsilon  = {10^{ - 6}},$ ${\alpha _{\min }} = {10^{ - 30}},$ ${\alpha _{\max }} = {10^{30}},$  
 ${\xi _0} = 1.07,$ ${\xi _1} = 5\times10^{ - 5}/3,$ $\xi _2=0.8,\;$ $\xi _3=5,\;$  $\sigma _{\min }=10^{-30}, \;\sigma _{\max }={10^{ 3}}$, $\delta  = {10^{ -4}},\;  {c_1} = {10^{ - 9}},$ ${c_2} = 10^{ - 7}$,  $c=0.99$ and
\begin{align*}
\alpha _0 = \left\{ \begin{array}{l}
 2\frac{{\left| {{f_0}} \right|}}{{\left\| {{g_0}} \right\|^2}},\;\;\;\;\;\; \quad \quad \quad \quad \quad \quad \quad \quad\;\;\;\;\;\;\;\; \;\;\;{\rm{if}}{\left\| {{x_0}} \right\|_\infty } < {10^{ - 30}}\,{\rm{and}}\,\left| {{f_0}} \right| \ge {10^{ - 30}},\\
1.0,\quad \quad \quad \quad \quad \quad \quad \quad \quad \quad \quad \quad \quad \;\;\;\;{\rm{if}}{\left\| {{x_0}} \right\|_\infty } < {10^{ - 30}}\,{\rm{and}}\,\left| {{f_0}} \right| < {10^{ - 30}},\\
\min \left\{ {1.0,\max \left\{ {\frac{{{{\left\| {{x_0}} \right\|}_\infty }}}{{{{\left\| {{g_0}} \right\|}_\infty }}},\frac{1}{{{{\left\| {{g_0}} \right\|}_\infty }}}} \right\}} \right\},\quad\ {\rm{if}}{\left\| {{x_0}} \right\|_\infty } \ge {10^{ - 30}}\,{\rm{and}}\,{\left\| {{g_0}} \right\|_\infty } \ge {10^7},\\
\min \left\{ {1.0,\frac{{{{\left\| {{x_0}} \right\|}_\infty }}}{{{{\left\| {{g_0}} \right\|}_\infty }}}} \right\},\quad \quad \quad \quad \quad \quad \quad \;\;\;{\rm{if}}{\left\| {{x_0}} \right\|_\infty } \ge {10^{ - 30}}\,{\rm{and}}\,{\left\| {{g_0}} \right\|_\infty } < {10^7}.
\end{array} \right.
\end{align*}

\begin{figure}[htp]
	\centering
	\begin{minipage}[t]{0.49 \linewidth}
		\includegraphics[scale=0.5]{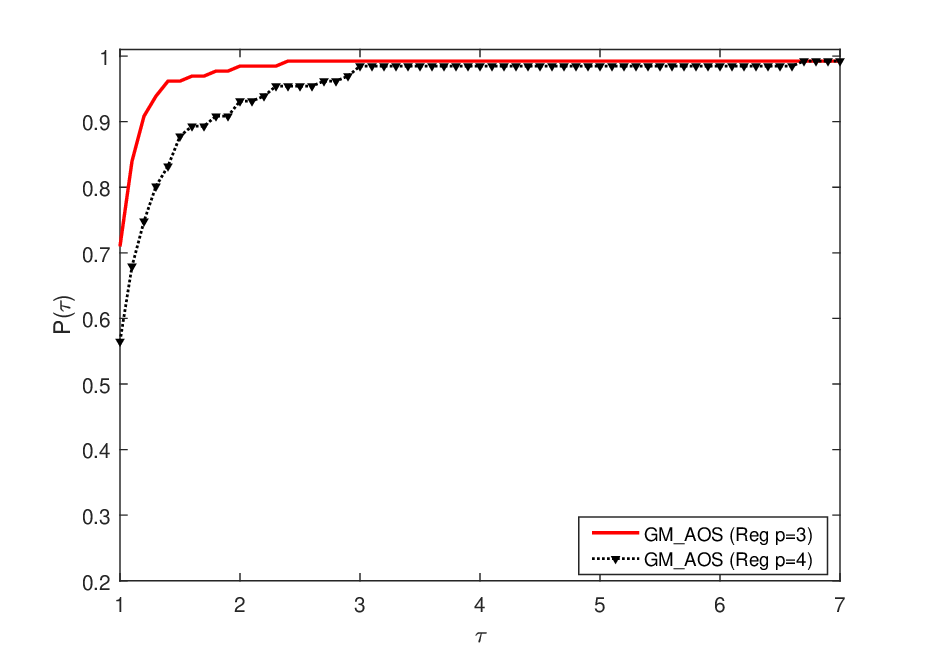}
		\caption{Performance profile based on $ N_{iter} $(CUTEr145)}\label{fig1}
	\end{minipage}	
	\begin{minipage}[t]{0.49 \linewidth}
		\centering
		\includegraphics[scale=0.5]{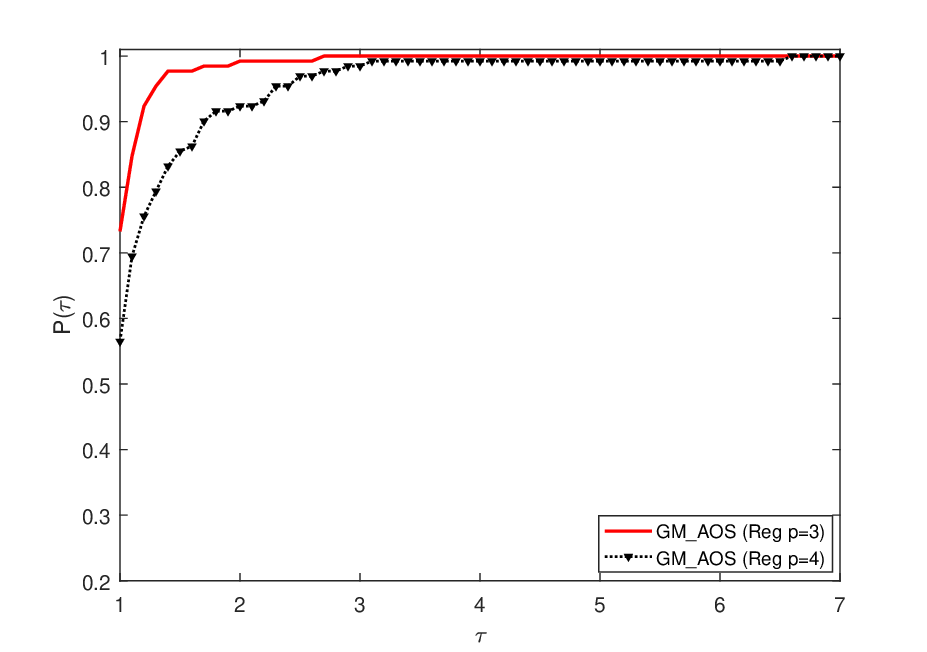}
		\caption{Performance profile based on $ N_f $(CUTEr145)}\label{fig2}
	\end{minipage}	
\end{figure}
\begin{figure}[htp]
	\centering
	\begin{minipage}[t]{0.49 \linewidth}
		\includegraphics[scale=0.5]{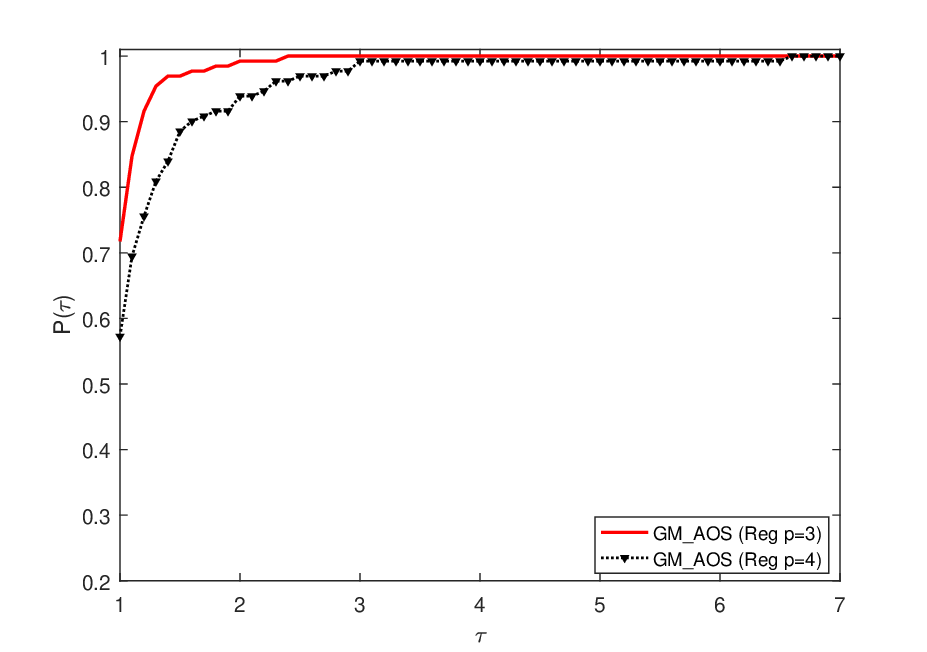}
		\caption{Performance profile based on $  N_g $(CUTEr145)}\label{fig3}
	\end{minipage}	
	\begin{minipage}[t]{0.49 \linewidth}
		\centering
		\includegraphics[scale=0.5]{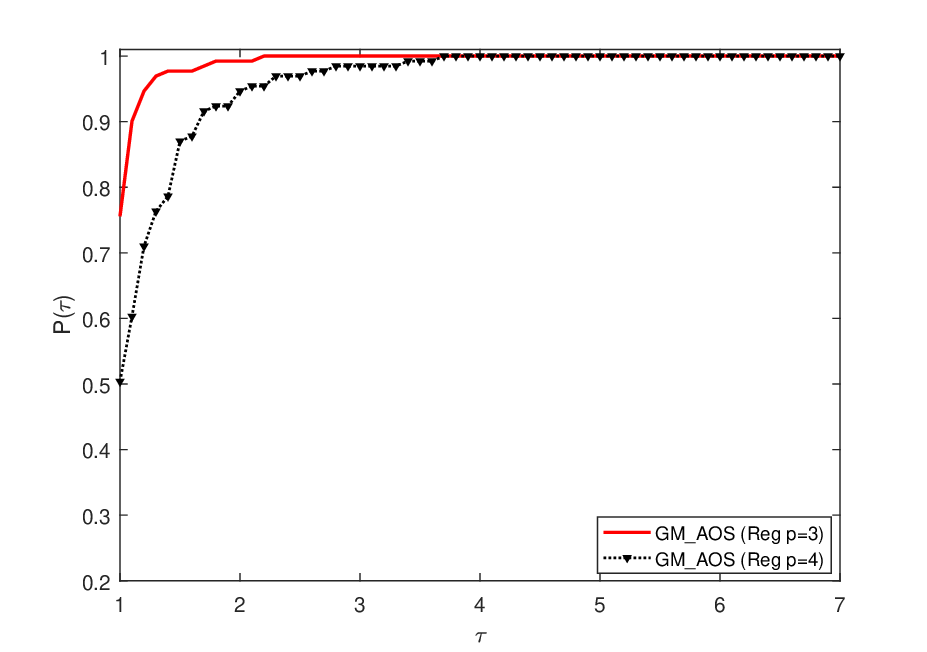}
		\caption{Performance profile based on $ T_{cpu} $(CUTEr145)}\label{fig4}
	\end{minipage}	
\end{figure}

\noindent  GM\_AOS (1.2)\cite{Liu2018Several} and the BB method used the same line search as that in  GM$ \_ $AOS (Reg). CGOPT(1.0),   CG$\_$DESCENT (5.0) and BBQ  used  all default settings of  parameters   but the   stopping conditions. Each test method is terminated if ${\left\| {{g_k}} \right\|_\infty } \le {10^{ - 6}} $ or the iterations exceeds 140,000.

The performance profiles introduced by Dolan and Mor\'e\cite{25} were used to display the performance of these methods.  In the following figures,  ``$ N_{iter} $'', ``$ N_f $'',
``$ N_g $'' and ``$ T_{cpu} $'' represent the number of
iterations,  the number of function evaluations, the number of gradient evaluations and CPU time (s), respectively.

We first compare GM$ \_ $AOS (Reg  p=3) with
GM$ \_ $AOS (Reg p=4) on the test set CUTEr145, and use the better one to  compare  with other test methods. As shown in Figs. \ref{fig1}-\ref{fig4}, we see that GM$ \_ $AOS (Reg p=3) performs better than GM$ \_ $AOS (Reg p=4) in term of   $ N_{iter} $,  $ N_f $,
 $ N_g $  and  $ T_{cpu} $. So we select GM$ \_ $AOS (Reg p=3)  to compare with other test methods in the following numerical experiments.

   \begin{figure}[htp]
   	\centering
   	\begin{minipage}[t]{0.49 \linewidth}
   		\includegraphics[scale=0.5]{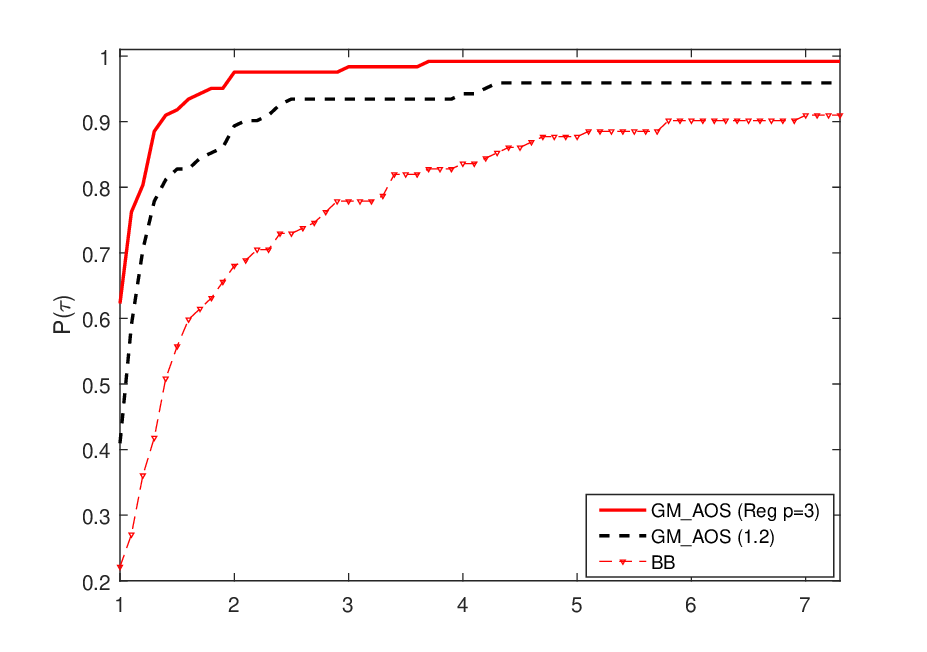}
   		\caption{Performance profile based on $ N_{iter}  $(CUTEr145)}\label{fig5}
   	\end{minipage}	
   	\begin{minipage}[t]{0.5 \linewidth}
   		\centering
   		\includegraphics[scale=0.5]{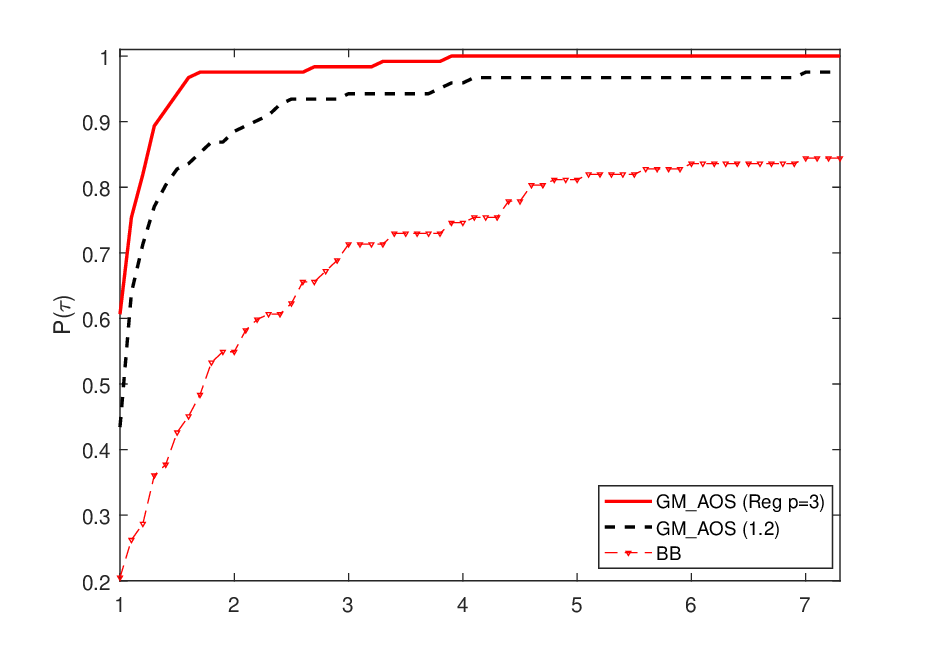}
   		\caption{Performance profile based on $ N_f  $(CUTEr145)}\label{fig6}
   	\end{minipage}	
   \end{figure}

      \begin{figure}[htp]
      	\centering
      	\begin{minipage}[t]{0.49 \linewidth}
      		\includegraphics[scale=0.5]{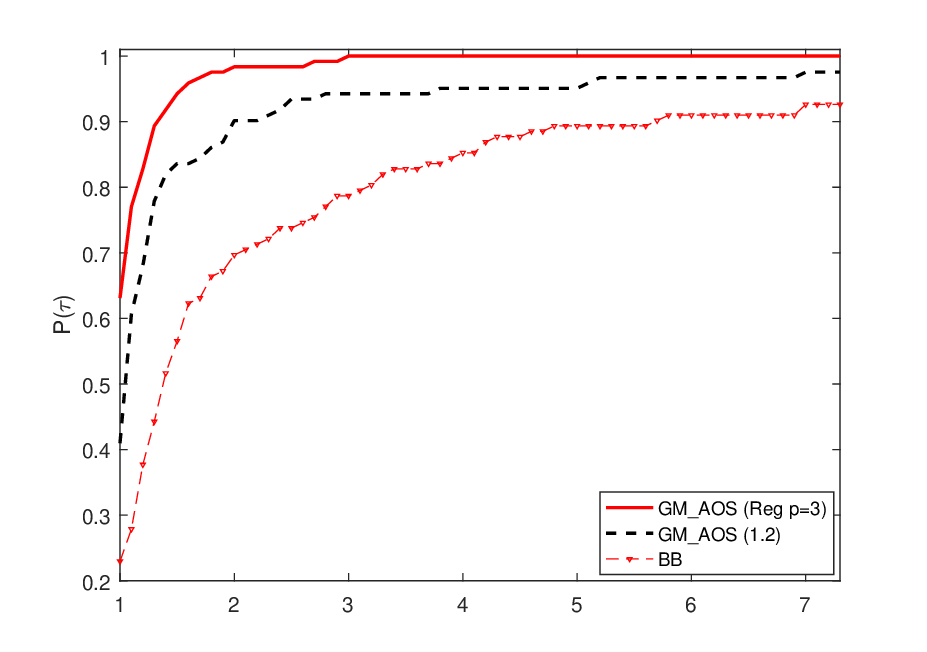}
      		\caption{Performance profile based on $   N_g $(CUTEr145).}\label{fig7}
      	\end{minipage}	
      	\begin{minipage}[t]{0.49 \linewidth}
      		\centering
      		\includegraphics[scale=0.45]{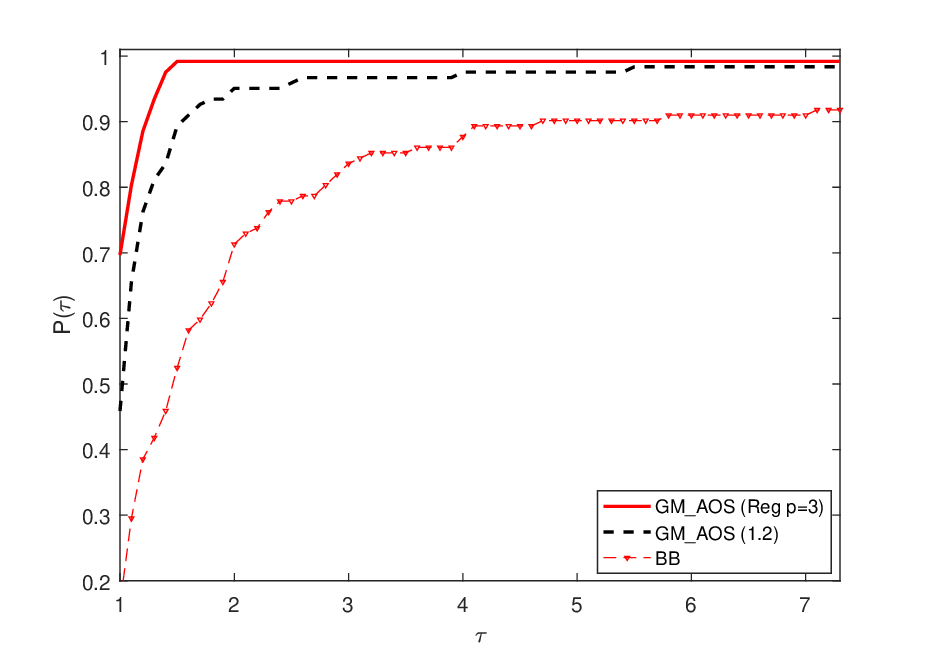}
      		\caption{Performance profile based on $ T_{cpu} $(CUTEr145).}\label{fig8}
      	\end{minipage}	
      \end{figure}

   The following numerical experiments are divided into  four groups.

  In the first group of the numerical experiments, we  compare the performance of GM$ \_ $AOS (Reg p=3) with that of GM$ \_ $AOS (1.2) \cite{Liu2018Several} and the BB method on the  test set  CUTEr145. Figs. \ref{fig5}-\ref{fig8} present the performance profiles  on the  test set  CUTEr145. As shown in Figs. \ref{fig5}-\ref{fig8}, we can observe that  GM$ \_ $AOS (Reg p=3) performs better than GM$ \_ $AOS (1.2) and is superior much to the BB method, and  GM$ \_ $AOS (1.2) outperforms   the BB method.   The first group of the numerical experiments indicates that the approximately optimal stepsize is   extremely  efficient.
  
  In the third group of the numerical experiments, we compare the numerical performance of GM$ \_ $AOS (Reg p=3) and  the BBQ method on the test set CUTEr145. In the numerical experiment, we do not compare the performance about the running time due to the fact that the BBQ method was implemented by Matlab code and GM$ \_ $AOS (Reg p=3) was implented by C code.  As shown in Fig. \ref{fig25}, \ref{fig26} and \ref{fig27}, we can observed that  GM$ \_ $AOS (Reg p=3) is superior to the  BBQ method for the test set CUTEr145 in term of the number of iteration, the number of funciton evaluation and the number of gradient evaluation, while the BBQ method has been regarded as the import advance for gradient method.
  
  \begin{figure}[htp]
  	\centering
  	\begin{minipage}[t]{0.49 \linewidth}
  		\includegraphics[scale=0.73]{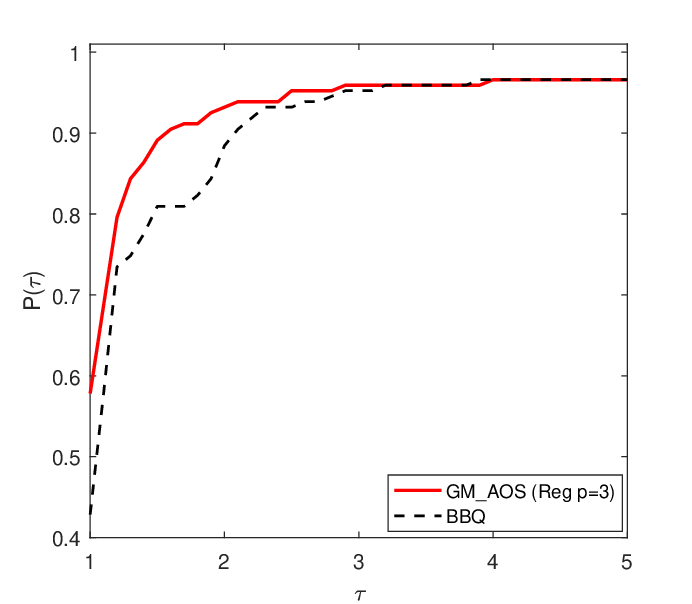}
  		\caption{Performance profile based on $ N_{iter} $}\label{fig25}
  	\end{minipage}	
  	\begin{minipage}[t]{0.49 \linewidth}
  		\centering
  		\includegraphics[scale=0.73]{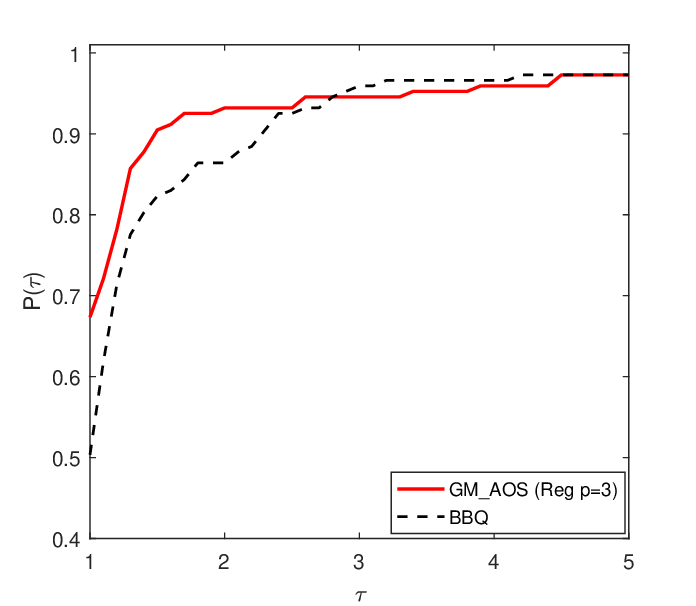}
  		\caption{Performance profile based on $ N_f $}\label{fig26}
  	\end{minipage}	
  \end{figure}
  \begin{figure}[htp]
  	\centering
  	\includegraphics[scale=0.8]{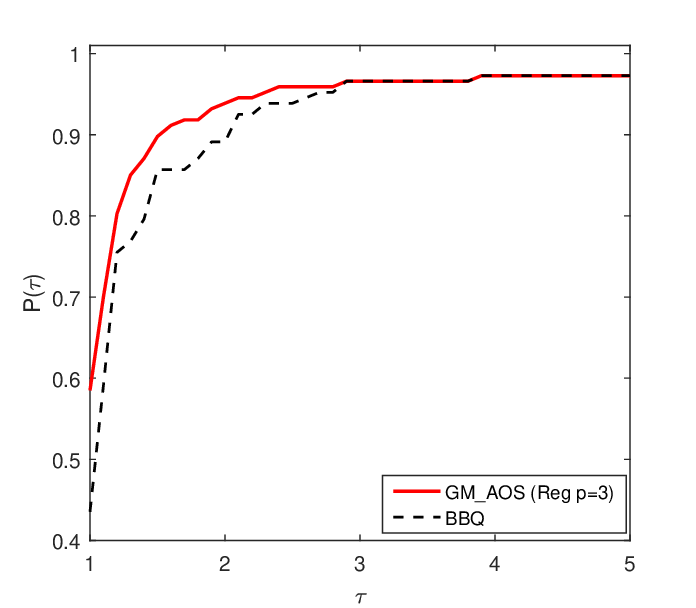}
  	\caption{Performance profile based on $  N_g $}\label{fig27}	
  	
  \end{figure}

  \begin{figure}[htp]
  	\centering
  	\begin{minipage}[t]{0.49 \linewidth}
  		\includegraphics[scale=0.5]{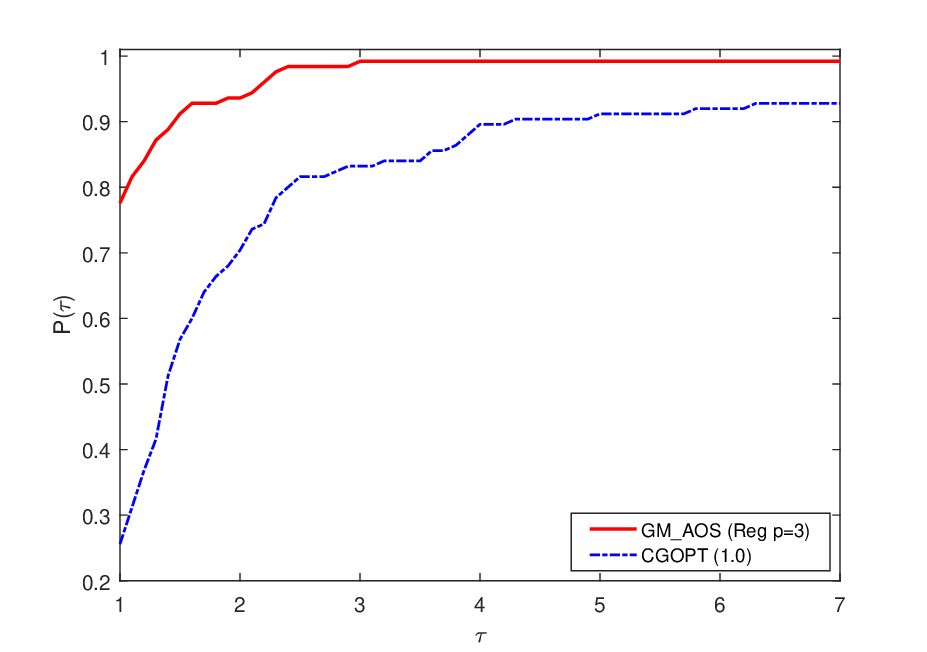}
  		\caption{Performance profile based on $ N_f  $(CUTEr145)}\label{fig9}
  	\end{minipage}	
  	\begin{minipage}[t]{0.5 \linewidth}
  		\centering
  		\includegraphics[scale=0.5]{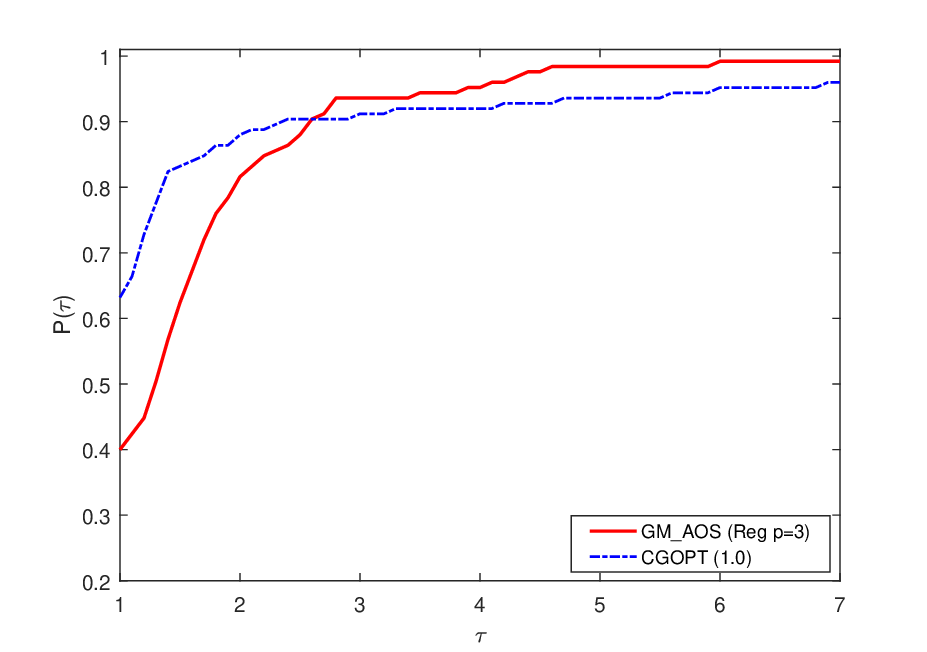}
  		\caption{Performance profile based on $ N_g  $(CUTEr145)}\label{fig10}
  	\end{minipage}	
  \end{figure}
  \begin{figure}[htp]
  	\centering
  	\begin{minipage}[t]{0.49 \linewidth}
  		\includegraphics[scale=0.49]{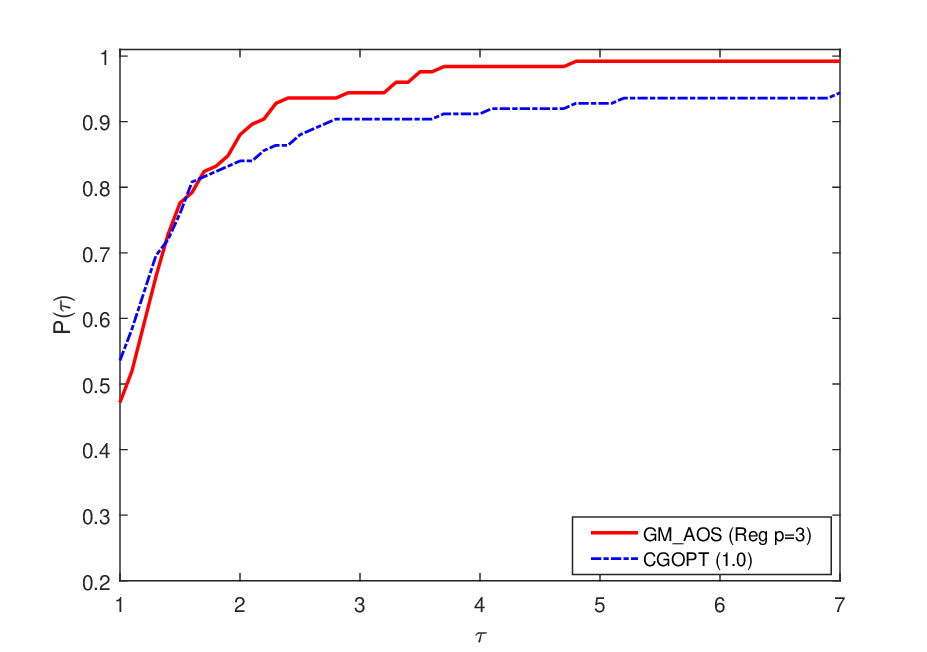}
  		\caption{Performance profile based on $N_f+3N_g $(CUTEr145).}\label{fig11}
  	\end{minipage}	
  	\begin{minipage}[t]{0.49 \linewidth}
  		\centering
  		\includegraphics[scale=0.49]{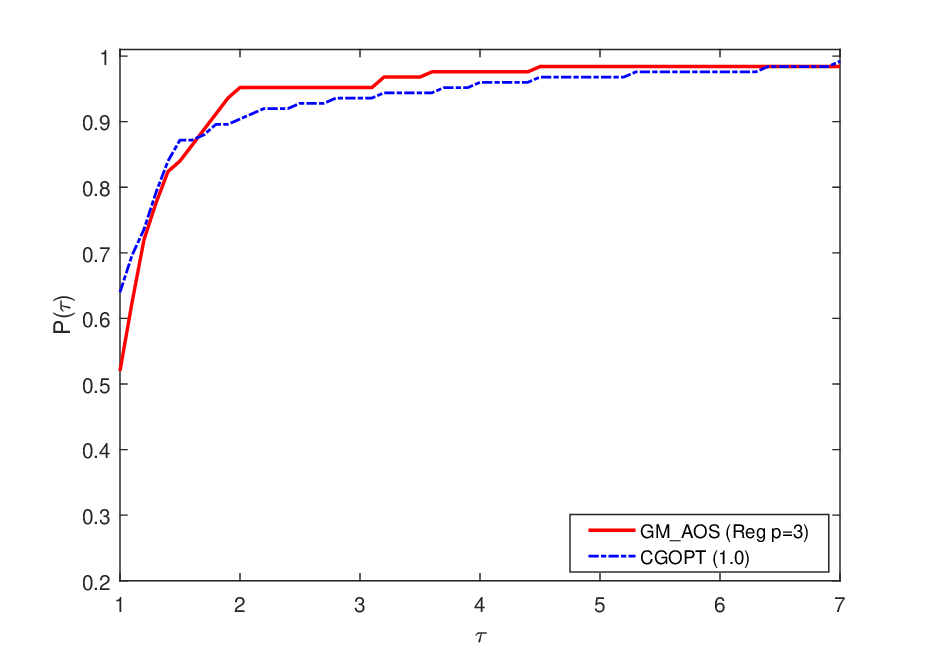}
  		\caption{Performance profile based on $ T_{cpu} $(CUTEr145).}\label{fig12}
  	\end{minipage}	
  \end{figure}

  \begin{figure}[htp]
  	\centering
  	\begin{minipage}[t]{0.49 \linewidth}
  		\includegraphics[scale=0.5]{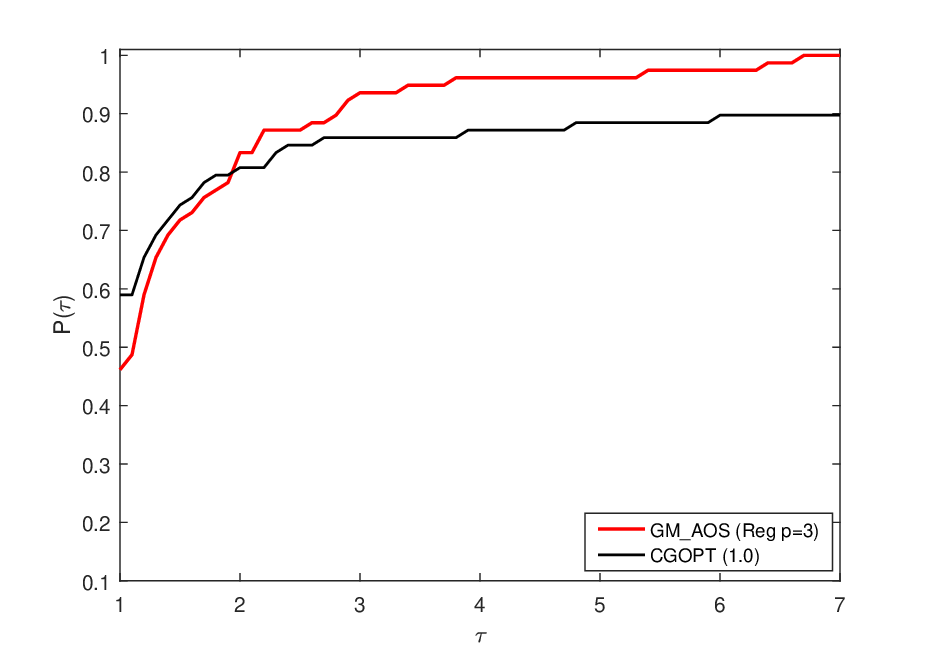}
  		\caption{Performance profile based on $ N_{iter}  $(Andr80)}\label{fig13}
  	\end{minipage}	
  	\begin{minipage}[t]{0.5 \linewidth}
  		\centering
  		\includegraphics[scale=0.5]{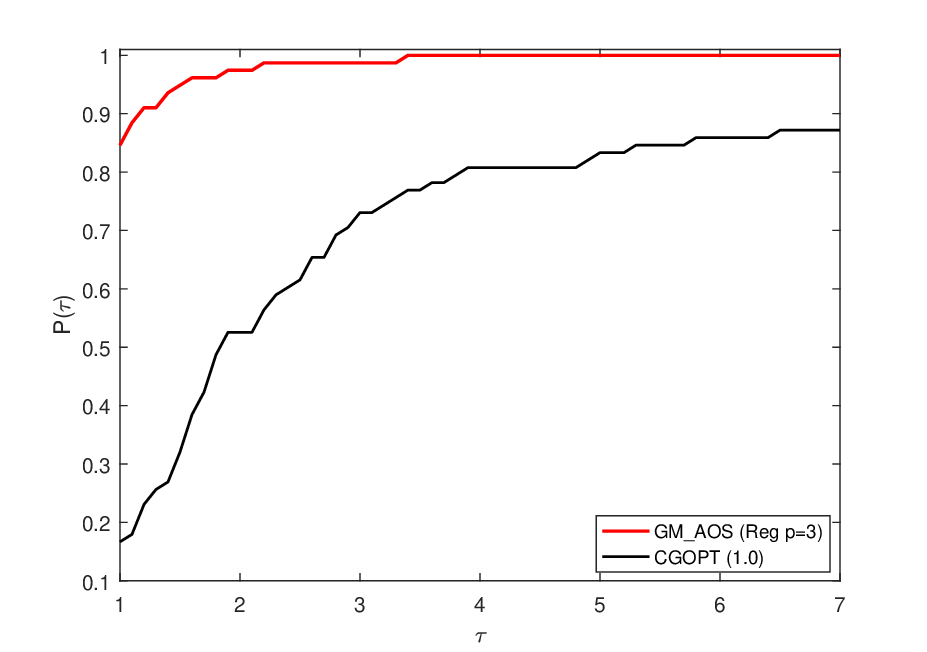}
  		\caption{Performance profile based on $ N_f  $(Andr80)}\label{fig14}
  	\end{minipage}	
  \end{figure}
  \begin{figure}[htp]
  	\centering
  	\begin{minipage}[t]{0.49 \linewidth}
  		\includegraphics[scale=0.5]{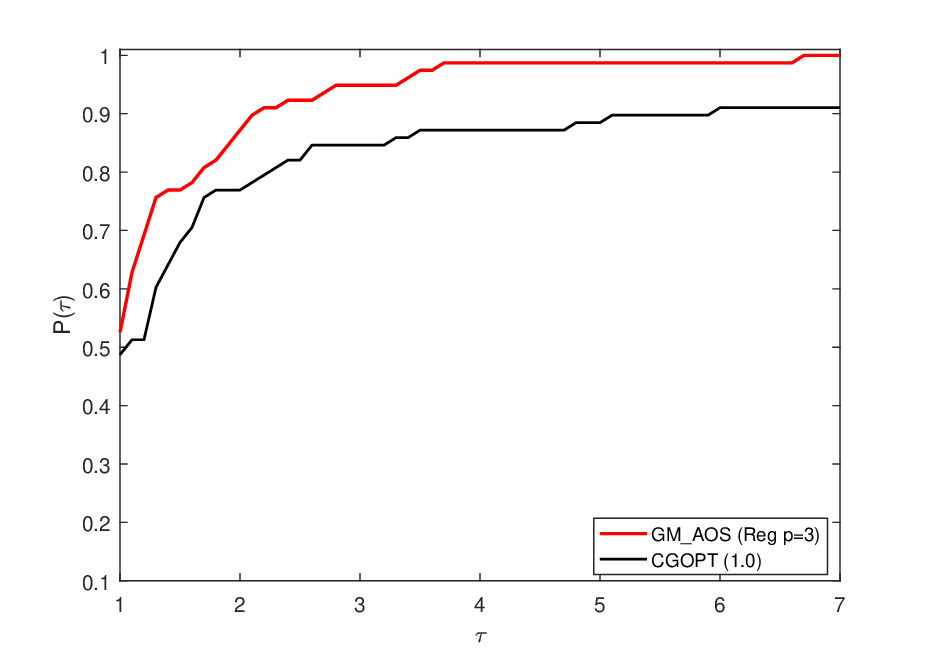}
  		\caption{Performance profile based on $  N_g $(Andr80).}\label{fig15}
  	\end{minipage}	
  	\begin{minipage}[t]{0.49 \linewidth}
  		\centering
  		\includegraphics[scale=0.55]{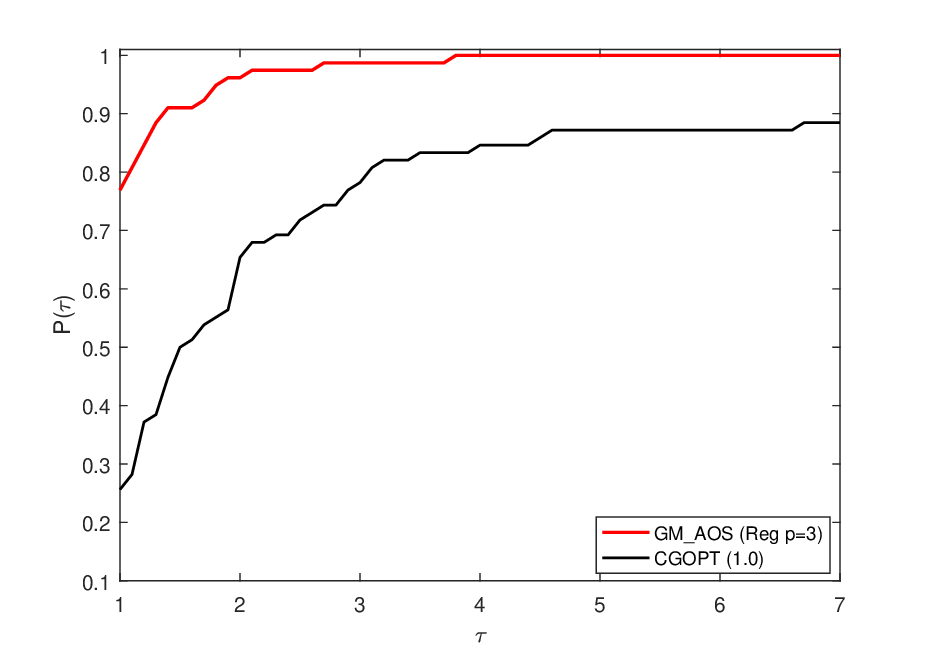}
  		\caption{Performance profile based on $ T_{cpu} $(Andr80).}\label{fig16}
  	\end{minipage}	
  \end{figure}

In the second group of the numerical experiments, we compare the performance of GM\_AOS (Reg p=3) with that of  CGOPT (1.0) on the two test sets  CUTEr145 and Andr80. Figs. \ref{fig9}-\ref{fig12} present the performance profiles on the  test set  CUTEr145.   As shown in Fig. \ref{fig9}, we see that GM\_AOS (Reg p=3) performs much better   CGOPT (1.0) in term of $N_f$, since GM\_AOS (Reg p=3) solves successfully about 79$\%$ test problems with the least function evaluations, while the percentage of CGOPT (1.0) is only about 38$\%$. Fig. \ref{fig10} indicates that GM\_AOS (Reg p=3)  is at a disadvantage over CGOPT (1.0) in term of $N_g$, and  Fig. \ref{fig11} shows that GM\_AOS (Reg p=3) outperforms slightly CGOPT (1.0) in term of $N_f +3N_g$ \cite{HagerLMCGDESCENT}. We can observe from Fig. \ref{fig12} that  GM\_AOS (Reg p=3) is as fast as   CGOPT (1.0). Figs. \ref{fig13}-\ref{fig16} present the performance profiles on the  test set  Andr80.  As shown in Figs. \ref{fig13}-\ref{fig16}, we observe that GM\_AOS (Reg p=3) illustrates
huge advantage over CGOPT (1.0) on the test set Andr80.  The second group of the numerical experiments indicates that GM\_AOS (Reg p=3)  is competitive to  CGOPT (1.0) on the   test set  CUTEr145, and has a significant improvement over  CGOPT (1.0) on the  test set  Andr80.

\begin{figure}[htp]
	\centering
	\begin{minipage}[t]{0.49 \linewidth}
		\includegraphics[scale=0.5]{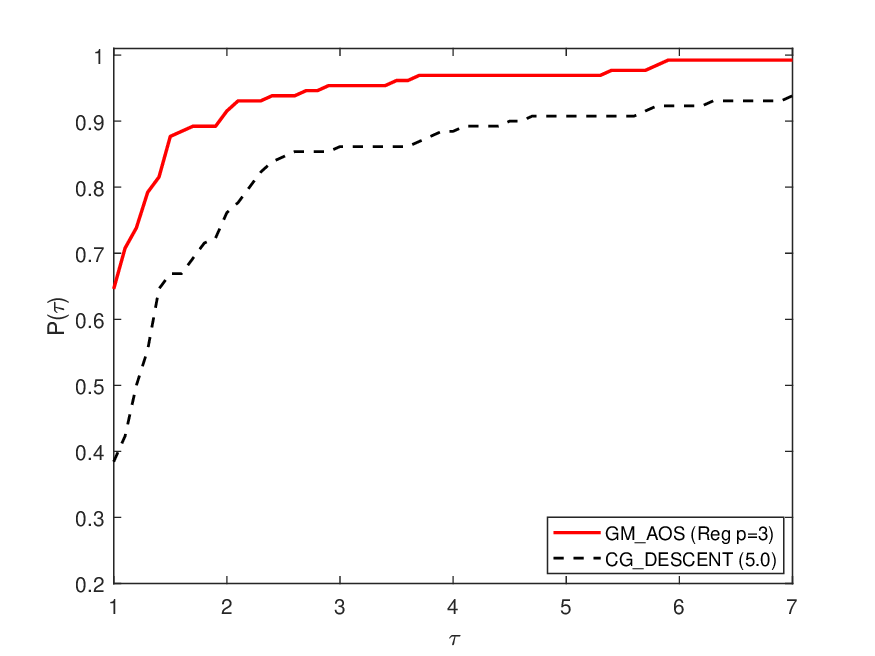}
		\caption{Performance profile based on $ N_{f} $(CUTEr145)}\label{fig17}
	\end{minipage}	
	\begin{minipage}[t]{0.49 \linewidth}
		\centering
		\includegraphics[scale=0.5]{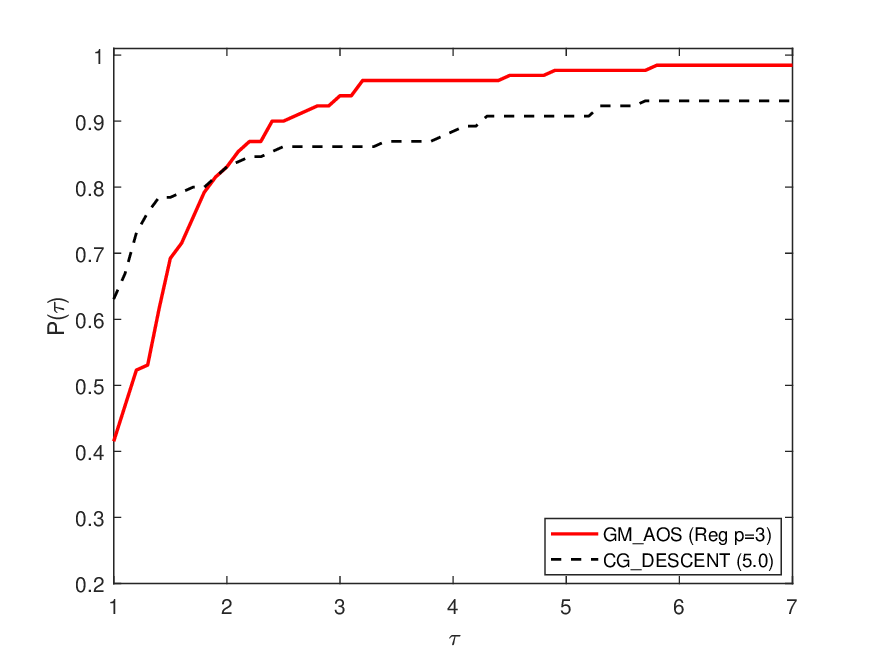}
		\caption{Performance profile based on $ N_g $(CUTEr145)}\label{fig18}
	\end{minipage}	
\end{figure}
\begin{figure}[htp]
	\centering
	\begin{minipage}[t]{0.49 \linewidth}
		\includegraphics[scale=0.5]{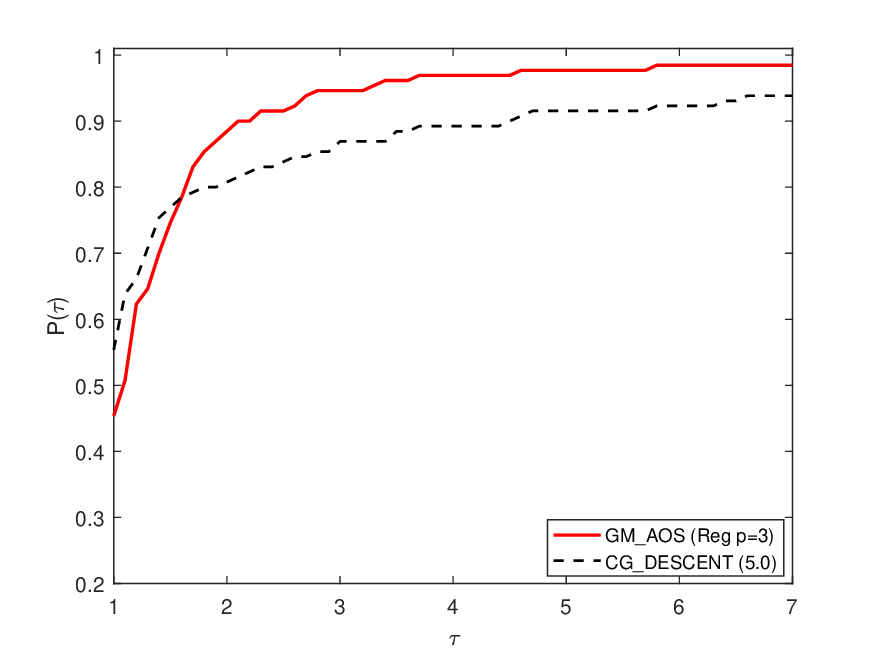}
		\caption{Performance profile based on $ N_{f}+3N_g $(CUTEr145)}\label{fig19}
	\end{minipage}	
	\begin{minipage}[t]{0.49 \linewidth}
		\centering
		\includegraphics[scale=0.5]{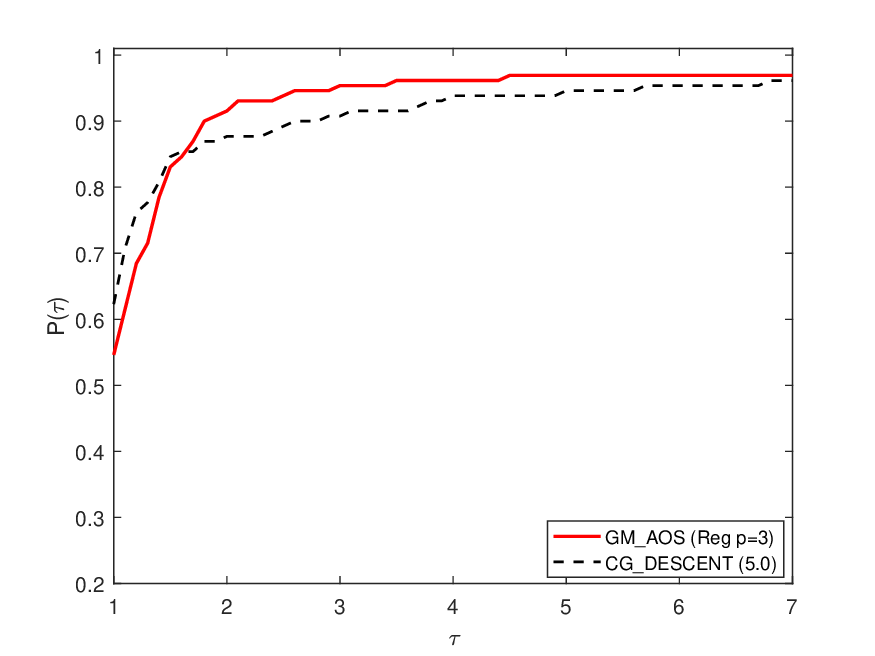}
		\caption{Performance profile based on $ T_{cpu} $(CUTEr145)}\label{fig20}
	\end{minipage}	
\end{figure}

\begin{figure}[htp]
	\centering
	\begin{minipage}[t]{0.49 \linewidth}
		\includegraphics[scale=0.5]{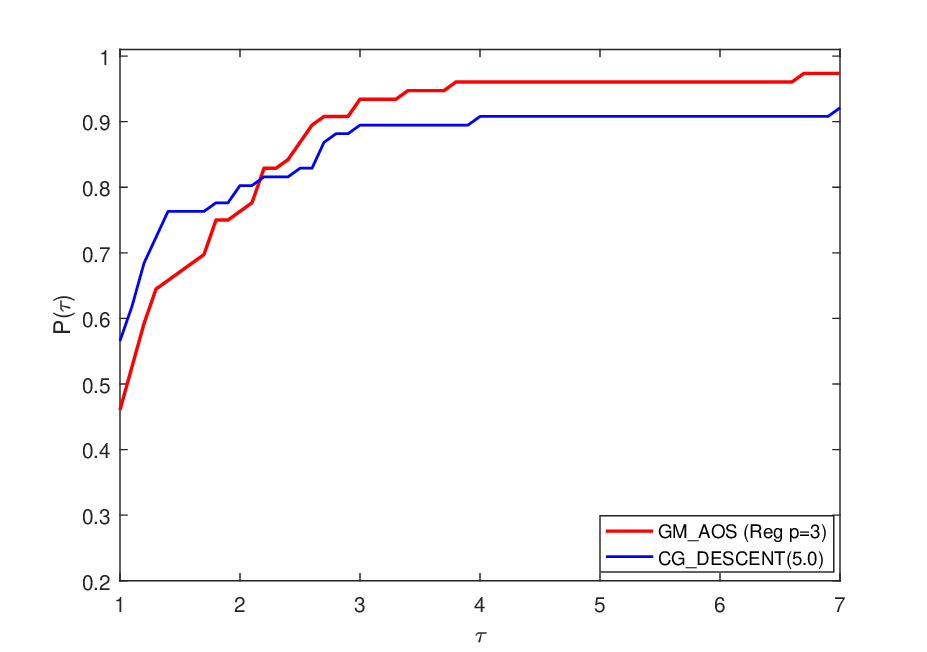}
		\caption{Performance profile based on $ N_{iter} $(Andr80)}\label{fig21}
	\end{minipage}	
	\begin{minipage}[t]{0.49 \linewidth}
		\centering
		\includegraphics[scale=0.5]{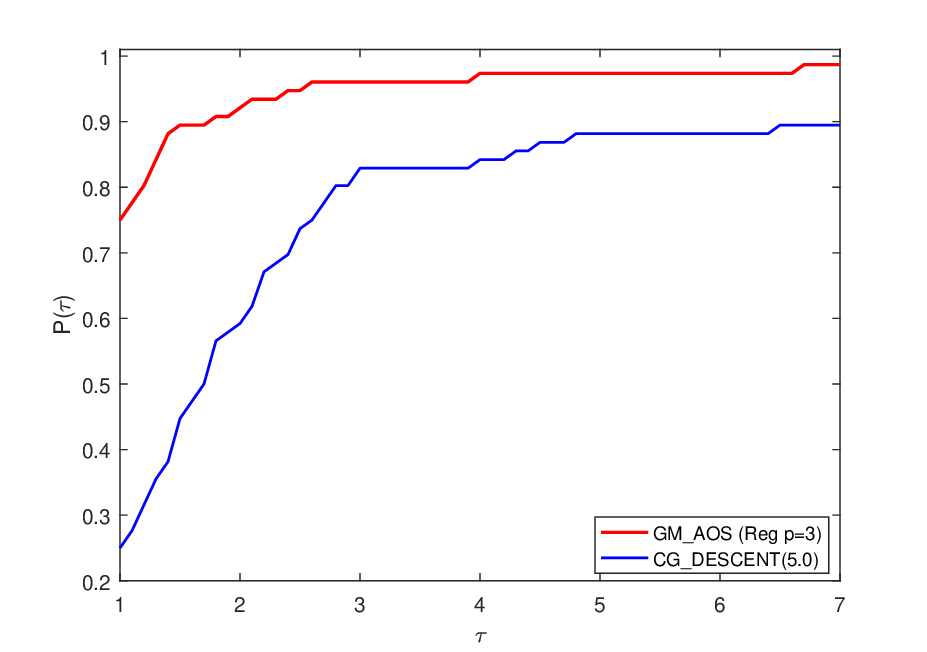}
		\caption{Performance profile based on $ N_f $(Andr80)}\label{fig22}
	\end{minipage}	
\end{figure}
\begin{figure}[htp]
	\centering
	\begin{minipage}[t]{0.49 \linewidth}
		\includegraphics[scale=0.5]{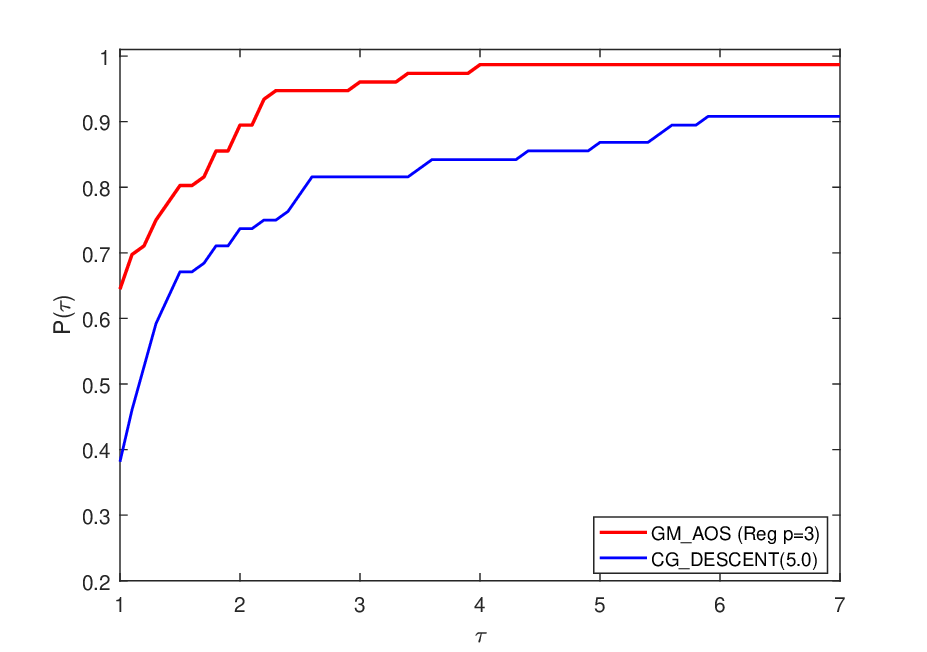}
		\caption{Performance profile based on $  N_g $(Andr80)}\label{fig23}
	\end{minipage}	
	\begin{minipage}[t]{0.49 \linewidth}
		\centering
		\includegraphics[scale=0.5]{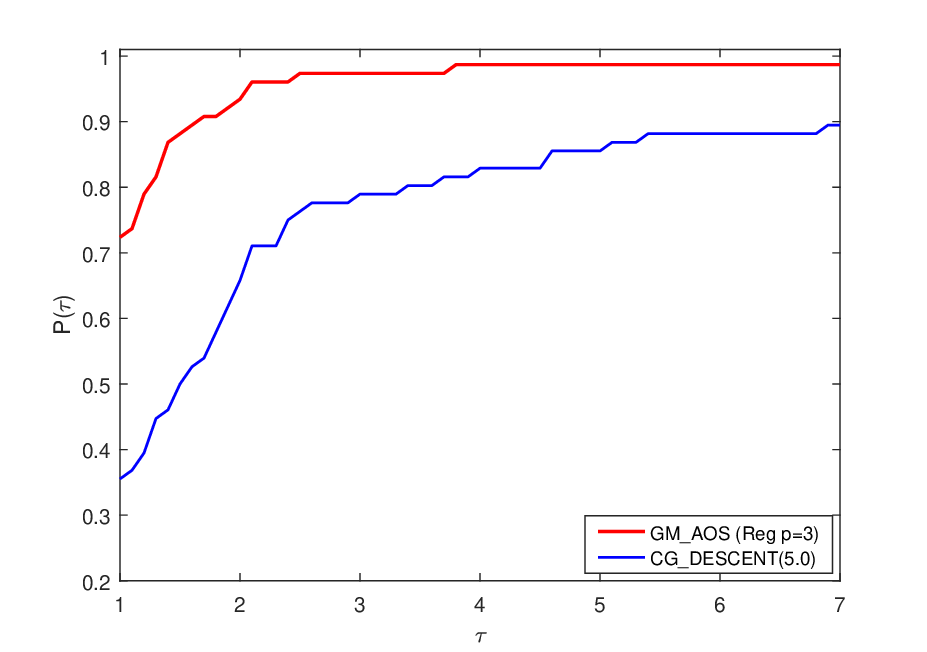}
		\caption{Performance profile based on $ T_{cpu} $(Andr80)}\label{fig24}
	\end{minipage}	
\end{figure}

In the fourth group of the numerical experiments,  we compare the performance of GM\_AOS (Reg p=3) with that of CG\_DESCENT (5.0) on the two test sets  CUTEr145 and Andr80.
Figs. \ref{fig17}-\ref{fig20} present the performance profiles  on the  test set  CUTEr145.  As shown in Fig. \ref{fig17}, we see that GM\_AOS (Reg p=3) performs better than CG\_DESCENT (5.0) in term of $N_f$,  since GM\_AOS (Reg p=3) solves successfully about 65$\%$ test problems with the least function evaluations, while the percentage of  CG\_DESCENT (5.0) is only about 39$\%$.   Fig. \ref{fig18} shows that GM\_AOS (Reg p=3) is at a disadvantage over than CG\_DESCENT (5.0) in term of $N_g$, and Fig. \ref{fig19} indicates  that GM\_AOS (Reg p=3) outperforms slightly CG\_DESCENT (5.0) in term of $N_f +3N_g$ \cite{HagerLMCGDESCENT}. We can observe from Fig. \ref{fig20} that GM\_AOS (Reg p=3)  is as fast as CG\_DESCENT (5.0).    Figs. \ref{fig21}-\ref{fig24} present the performance profiles  on the  test set  Andr80.  As shown in Figs. \ref{fig21}-\ref{fig23}, we see that GM\_AOS (Reg p=3) is at a little disadvantage over  CG\_DESCENT (5.0) in term of $N_{iter}$, and has a   significant performance boost over  CG\_DESCENT (5.0) in term of $N_f$  and  $N_g$. We also can see that  GM\_AOS (Reg p=3) is faster much than CG\_DESCENT (5.0).  The third group of the numerical experiments indicates that GM\_AOS (Reg p=3)  is competitive to  CG\_DESCENT (5.0) on the   test set  CUTEr145, and has a significant improvement over CG\_DESCENT (5.0) on the  test set  Andr80. \newline   

\begin{center}
	\small{\small  Table 1.~The number of test problems   }\\[2.0mm]
	\renewcommand\arraystretch{1.25}
	\begin{longtable}{{lccc cccc c c }}
		\toprule	\label{table1}
		Method &	$N_{\text{linsear}}=0$	&	$N_{\text{linsear}}\le1$	&	$N_{\text{linsear}}\le2$ &	$N_{\text{linsear}}\le 3$ 	& total problems	\\
		\midrule		
		BB	&	 41 	&	46	& 48	&50	&	145(CUTEr145) & 	\\
		GM\_AOS (Reg p=3)	&	\textbf{68}	&	81	& 85	& \textbf{90} 	&	 	145(CUTEr145) & 	\\
		GM\_AOS (Reg p=4)	&	\textbf{51}	&	57	& 60	&\textbf{62}	&	 	80(Andr80) & 	 \\	
		\bottomrule	
	\end{longtable}
\end{center}

As for the reasons for the suprising numerical effect of GM\_AOS (Reg p=3), we think  that they lie  in two aspects: (i)The approximately optimal stepsize is generated by the approximation models including regularization models and quadratic models   at the current iterate $x_k$,  which implies that it  is  incorporated  

\noindent properly into more second order or  high order information of the objective function.  (ii)The approximately optimal stepsize can   readily satisfy  Zhang-Hager line search directly in most cases   compared to other stepsizes in gradient method, which implies that it   require less much function evaluations and thus save  much computational cost. This  can be observed in Figs. \ref{fig6}, \ref{fig9}, \ref{fig14}, \ref{fig17} and  \ref{fig22}. More     results can be seen in Table \ref{table1}.  In Table \ref{table1}, $N_{\text{linsear}}$ represents  the times that the stepsize is  updated by \eqref{eq:ZhangHageLS4}   during all iterations of solving a test problem, namely, the times that   Zhang-Hager line search   is invoked  during \textbf{all iterations} of solving a test problem.
 $N_{\text{linsear}}=0$ indicates the initial stepsize (approximately optimal stepsize or BB stepsize) satisfies Zhang-Hager line search \eqref{eq:ZhangHageLS1}  directly at all iterations and thus  \textbf{Zhang-Hager line search }  
    	\textbf{is not invoked at all}. As shown in Table \ref{table1},  we can see that   there are 68  (out of 145) problems   for which  Zhang-Hager line search is not invoked at all during the solving process, while the number for the BB method is only 41, and there are 90   (out of 145)   problems  for each of which   the times that Zhang-Hager line search is invoked is  less than  or equal  to 3, while the number  for the BB method is only 50. We also can see that there are  51 (out of 80) problems  for which  Zhang-Hager line search is not invoked at all during the solving process. Table \ref{table1} indicates that the  approximately optimal stepsize in GM\_AOS (Reg p=3) is easier much to meet Zhang-Hager line search \eqref{eq:ZhangHageLS1} directly.  


\section{Conclusion and discussion}

In this paper, we present two     efficient gradient methods with approximately optimal stepsize  for unconstrained optimization. In the proposed method, some approximation models including   regularization models and quadratic models are exploited carefully to derive    approximately optimal stepsize. The convergence of the proposed methods is analyzed.  Extensive numerical results indicates that the proposed method GM$\_$AOS (Reg p=3) is superior to the  BBQ method   and other efficient gradient methods, and   is competitive to two  quite efficient and well-known conjugate gradient  software packages CG$ \_ $DESCENT (5.0)  and CGOPT (1.0) on the 145 test problems in the CUTEr library, has significant improvement over  CG$ \_ $DESCENT (5.0) and  CGOPT (1.0) on the  80 test problems collected by Andrei. As for the reason that GM$\_$AOS (Reg p=3) has so important improvement over  CG$ \_ $DESCENT (5.0) and  CGOPT (1.0) on Andr80 and is only competitive to  CG$ \_ $DESCENT (5.0) and  CGOPT (1.0) on CUTEr145, I think that it lies mainly in that   most     test problems in CUTEr145 is relatively difficult to solve compared to  the  test problems in Andr80. It seems that one can draw the following conclusion: Gradient methods with approximately optimal stepsize are sufficient  for those test problems that are not very ill-conditioned.

 Given that the facts:  (i)the search direction $-g_k$     has low storage; (ii)the approximately optimal stepsize can be easily computed;  (iii)the  nonmonotone Armijo  line search used   can be easily implemented; (iv)the numerical effect is   surprisingly nice, the  gradient methods with approximately optimal stepsizes can  become   strong candidates  for  large scale  unconstrained optimization and has potential in constrained optimization and some fields such as  machine learning.
 
 Though gradient methods with approximately optimal stepsize is surprisingly efficient, there are still some questions under investigation:
 
(i)Does gradient method  with approximately optimal stepsize based on quadratic approximation model  \eqref{eq:Quadmodel4}   possess Q-linear convergence for convex quadratic minimization? If yes, what conditions should be imposed on the distance $\left\| B_k - A \right\|  $ ? Here $A$ is the Hessian matrix for strictly   convex quadratic function. 

(ii)It will be an interesting   research for combining approximately optimal stepsize wih Cauchy stepsize   in  convex quadratic minimization.   How should one combine approximately optimal stepsize with Cauchy stepsize    for obtaining better convergence rate in   convex quadratic minimization?

(iii)Can the type of gradient method  with approximately optimal stepsize  possess local  R-linear convergence or better convergence rate   when  it is  applied to   general unconstrained optimization?

(iv)There are still a large room for improving numerical performance of gradient methods with approximately optimal stepsizes by exploiting other adaptive appriximation models based on the properties of the objective function.

\begin{acknowledgements} 
  We  would like to   thank Professors Hager W. W. and Zhang, H. C. for their C code of CG$ \_ $DESCENT, thank Professor Dai Y. H. and Dr. Kou C. X. for their C code of CGOPT (1.0). This research is supported by National Science Foundation of China (No. 11901561) and Guangxi Natural Science Foundation (2018GXNSFBA281180).
\end{acknowledgements}

\end{document}